\documentclass{amsart}

\usepackage{amsfonts}
\usepackage{amssymb}
\usepackage{enumerate}
\usepackage{amsmath}
\usepackage[T1]{fontenc}
\usepackage{graphicx}
\usepackage{lmodern}
\usepackage[unicode]{hyperref}

\newtheorem{ttt}{Theorem}[section]
\newtheorem{llll}[ttt]{Lemma}
\newtheorem{ccc}[ttt]{Claim}
\newtheorem{eee}[ttt]{Example}
\newtheorem{fff}[ttt]{Fact}
\newtheorem{rrr}[ttt]{Remark}
\newtheorem{sss}[ttt]{Statement}
\newtheorem{ddd}[ttt]{Definition}
\newtheorem{qqq}[ttt]{Question}
\newtheorem{cccc}[ttt]{Corollary}
\newtheorem{nnn}[ttt]{Notation}
\newtheorem{ppp}[ttt]{Problem}
\newtheorem{pppp}[ttt]{Proposition}
\newtheorem{ccccc}[ttt]{Conjecture}

\newcommand{\beq}{\begin{equation} }

\newcommand{\bt}{\begin{ttt}}
\newcommand{\bl}{\begin{llll}}
\newcommand{\bc}{\begin{ccc}}
\newcommand{\bex}{\begin{eee}}
\newcommand{\bfa}{\begin{fff}}
\newcommand{\br}{\begin{rrr}\upshape}
\newcommand{\bst}{\begin{sss}}
\newcommand{\bd}{\begin{ddd}\upshape}
\newcommand{\bdd}{\begin{ddd}\upshape}

\newcommand{\bq}{\begin{qqq}}
\newcommand{\bnn}{\begin{nnn}}
\newcommand{\bpr}{\begin{ppp}}
\newcommand{\bprop}{\begin{pppp}}
\newcommand{\bcor}{\begin{cccc}}
\newcommand{\bcon}{\begin{ccccc}}

\newcommand{\eeq}{\end{equation}}
\newcommand{\et}{\end{ttt}}
\newcommand{\el}{\end{llll}}
\newcommand{\ec}{\end{ccc}}
\newcommand{\eex}{\end{eee}}
\newcommand{\efa}{\end{fff}}
\newcommand{\er}{\end{rrr}}
\newcommand{\est}{\end{sss}}
\newcommand{\ed}{\end{ddd}}
\newcommand{\eq}{\end{qqq}}
\newcommand{\ecor}{\end{cccc}}
\newcommand{\econ}{\end{ccccc}}
\newcommand{\enn}{\end{nnn}}
\newcommand{\epr}{\end{ppp}}
\newcommand{\eprop}{\end{pppp}}

\newcommand{\bp}{\noindent\textbf{Proof. }}
\newcommand{\ep}{\hspace{\stretch{1}}$\square$\medskip}

\newcommand{\lab}[1]{\label{#1}}

\newcommand{\ZZ}{\mathbb{Z}}

\newcommand{\RR}{\mathbb{R}}

\newcommand{\al}{\alpha}

\newcommand{\om}{\omega}
\newcommand{\si}{\sigma}
\newcommand{\ka}{\kappa}

\newcommand{\iA}{\mathcal{A}}
\newcommand{\iB}{\mathcal{B}}

\newcommand{\iF}{\mathcal{F}}

\newcommand{\iI}{\mathcal{I}}

\newcommand{\iM}{\mathcal{M}}

\newcommand{\iP}{\mathcal{P}}

\newcommand{\iN}{\mathcal{N}}

\newcommand{\beeq}{\begin{equation}}
\newcommand{\eeeq}{\end{equation}}
\def\su{\subset}

\newcommand{\non}{{\rm non}}
\newcommand{\cof}{{\rm cof}}
\newcommand{\add}{{\rm add}}
\newcommand{\cov}{{\rm cov}}
\newcommand{\cl}{{\rm cl}}
\newcommand{\dom}{{\rm dom}}
\newcommand{\ran}{{\rm ran}}

\newcommand{\diam}{\mathrm{diam}}

\newcommand{\HNG}{\mathcal{HN}(G)}
\newcommand{\HN}{\mathcal{HN}}
\newcommand{\HMG}{\mathcal{HM}(G)}
\newcommand{\HM}{\mathcal{HM}}

\newcommand{\tieconcat}{%
	\mathbin{\mathpalette\dotieconcat\relax}%
}
\newcommand{\dotieconcat}[2]{
	\text{\raisebox{.8ex}{$\smallfrown$}}%
}

\newcommand\restrict[1]{\raisebox{-.4ex}{$|$}_{#1}}
\numberwithin{equation}{section}

\begin{document}

\author[M\'arton Elekes]{M\'arton Elekes$^\ast$}
\address{Alfr\'ed R\'enyi Institute of Mathematics, Hungarian Academy of Sciences,
PO Box 127, 1364 Budapest, Hungary and E\"otv\"os Lor\'and
University, Institute of Mathematics, P\'azm\'any P\'eter s. 1/c,
1117 Budapest, Hungary}
\email{elekes.marton@renyi.mta.hu}
\urladdr{http://www.renyi.hu/$\sim$emarci}
\thanks{The authors were supported by the National Research, Development and Innovation Office
-- NKFIH, grants no. 104178, 124749 and 129211. The first author was also supported by the National Research, Development and Innovation Office
-- NKFIH, grant no. 113047.
\parbox[c]{12cm}{\scshape Supported by the \'UNKP-18-3 New National Excellence Program of the Ministry of Human Capacities.}
}

\author[M\'ark Po\'or]{M\'ark Po\'or$^\dag$}
\address{E\"otv\"os Lor\'and
University, Institute of Mathematics, P\'azm\'any P\'eter s. 1/c,
1117 Budapest, Hungary}
\email{sokmark@caesar.elte.hu}

 	


\title{Cardinal invariants of Haar null and Haar meager sets}

\begin{abstract}
A subset $X$ of a Polish group $G$ is \emph{Haar null} if there exists a Borel probability measure $\mu$ and a Borel set $B$ containing $X$ such that $\mu(gBh)=0$ for every $g,h \in G$. A set $X$ is \emph{Haar meager} if there exists a compact metric space $K$, a continuous function $f : K \to G$ and a Borel set $B$ containing $X$ such that $f^{-1}(gBh)$ is meager in $K$ for every $g,h \in G$. We calculate (in $ZFC$) the four cardinal invariants ($\rm add$, $\rm cov$, $\rm non$, $\rm cof$) of these two $\sigma$-ideals for the simplest non-locally compact Polish group, namely in the case $G = \mathbb{Z}^\omega$. In fact, most results work for separable Banach spaces as well, and many results work for Polish groups admitting a two-sided invariant metric. This answers a question of the first named author and Z. Vidny\'anszky.
\end{abstract}

\maketitle

\section{Introduction}

Small sets play a fundamental role in many branches of mathematics. Perhaps the most important example is the family of nullsets of a natural invariant measure. Such a natural measure is the Lebesgue measure on $\mathbb{R}^d$, or more generally the Haar measure on a locally compact group. However, on larger groups such as $C[0,1]$ or $S_\infty$ there is no such measure (a Polish group carries a Haar measure iff it is locally compact, see e.g. \cite{EN}). Therefore J. P. R. Christensen \cite{Ch} introduced the following notion.

\bd
A subset $X$ of a Polish group $G$ is \emph{Haar null} if there exists a Borel probability measure $\mu$ and a Borel set $B$ containing $X$ such that $\mu(gBh)=0$ for every $g,h \in G$.

The family of Haar null sets is denoted by $\HNG$ or simply $\HN$.
\ed

Christensen proved that these sets form a proper $\sigma$-ideal which coincides with the family of sets of Haar measure zero in the locally compact case. This notion turned out to be very useful in various branches of mathematics, see e.g. the survey paper \cite{EN}.

\bigskip
Actually, the right dual notion to Haar null sets is not the meager sets, hence U. B. Darji \cite{Da} introduced the following notion.

\bd
A subset $X$ of a Polish group $G$ is \emph{Haar meager} if there exists a compact metric space $K$, a continuous function $f : K \to G$ and a Borel set $B$ containing $X$ such that $f^{-1}(gBh)$ is meager in $K$ for every $g,h \in G$.

The family of Haar meager sets is denoted by $\HMG$ or simply $\HM$.
\ed  

Analogously to the Haar null case, Darji proved that these sets form a proper $\sigma$-ideal which coincides with the family of meager sets in the locally compact case.
For more information see e.g. the survey paper \cite{EN}.

\bigskip
When investigating a notion of smallness, a fundamental concept is that of cardinal invariants. The four most notable ones are the following.

\bd
Let $\iI$ be a $\si$-ideal on a set $X$. Define
\[
\begin{array}{lll}
\add(\iI)& =& \min\{|\iA| : \iA \su \iI, \bigcup \iA \notin \iI\},\\
\cov(\iI)& =& \min\{|\iA| : \iA \su \iI, \bigcup \iA = X\},\\
\non(\iI)& =& \min\{|H| : H \su X, H \notin \iI\},\\
\cof(\iI)& =& \min\{|\iA| : \iA \su \iI, \ \forall I \in \iI \ \exists A \in \iA, I \su A\}.
\end{array}
\]
These invariants are called the \emph{additivity, covering number, uniformity, and cofinality of $\iI$}, respectively.
\ed 
For more information on cardinal invariants see e.g. the monograph \cite{BJ}.

\bigskip
The goal of this paper is to determine these cardinal invariants of $\HN$ and $\HM$. The case of $\HN$ was asked in \cite[Question 5.7]{EV}. Before we proceed, let us describe the most important results on this topic so far. First, note that if $G$ is locally compact then the Haar null sets agree with the sets of Haar measure zero and Haar meager sets agree with the meager sets \cite{EN}. Moreover, it is also well-known that the four invariants of the measure zero sets do not depend on the underlying measure space as long as it is a Polish space equipped with a continuous $\si$-finite Borel measure, e.g. a locally compact non-discrete Polish group equipped with the (left) Haar measure. Similarly, the four invariants of the meager sets do not depend on the underlying space as long as it is a Polish space without isolated points, e.g. a non-discrete Polish group \cite{BJ}. Therefore, if $G$ is locally compact and non-discrete then the four invariants of $\HNG$ agree with the respective invariants of $\iN$ (the family of Lebesgue nullsets of $\RR$), and the four invariants of $\HMG$ agree with the respective invariants of $\iM$ (the family of meager subsets  of $\RR$).

\bigskip
Recall that a set is \emph{universally measurable} if it is measurable with respect to the completion of every Borel probability measure.
\bd
A subset $X$ of a Polish group $G$ is \emph{generalized Haar null} if there exists a (completed) Borel probability measure $\mu$ and a universally measurable set $B$ containing $X$ such that $\mu(gBh)=0$ for every $g,h \in G$.

The family of generalized Haar null sets is denoted by $\HN_{gen}(G)$ or simply $\HN_{gen}$.
\ed

The following results were proved by T. Banakh \cite{Ba}. For the definition of the so called bounding number $\mathfrak{b}$ and dominating number $\mathfrak{d}$ see e.g. \cite{BJ}.

\bt
\lab{t:gen}
\[
\begin{array}{lll}
\add(\HN_{gen}(\ZZ^\om))& = & \add(\iN),\\
\cov(\HN_{gen}(\ZZ^\om))& = & \min\{\mathfrak{b},  \cov(\iN)\},\\
\non(\HN_{gen}(\ZZ^\om))& = & \max\{ \mathfrak{d, \non(\iN)}\},\\
\smallskip 
\hspace{-90pt}\textrm{and under Martin's Axiom}&&\\
\cof(\HN_{gen}(\ZZ^\om))& > & \mathfrak{c}.
\end{array}
\]
\et

This last statement is really peculiar, since all the usual cardinal invariant are at most the continuum.

\bigskip
Now we turn to the results concerning $\HN$. The following theorem will answer \cite[Question 5.7]{EV}, and will show a surprising contrast to the above results of Banakh.

For the sake of completeness we list the values of all four invariants, but note that additivity was already calculated by the first named author and Z. Vidny\'anszky \cite{EV}.

\bt
\[
\begin{array}{lll}
\add(\HN(\ZZ^\om))& = & \om_1,\\
\cov(\HN(\ZZ^\om))& = & \min\{\mathfrak{b},  \cov(\iN)\},\\
\non(\HN(\ZZ^\om))& = & \max\{ \mathfrak{d, \non(\iN)}\},\\
\cof(\HN(\ZZ^\om))& = &\mathfrak{c}.
\end{array}
\]
\et

\br
In fact, additivity and cofinality works for all non-locally compact Polish groups admitting a two-sided invariant metric.
\er

\bigskip
Now we turn to the case of $\HM$.

Again, we also list additivity here, which was already calculated by M. Dole\v{z}al an V. Vlas\'ak \cite{DV}. Let $\iM$ denote the family of meager subsets of $\RR$.

\bt
\[
\begin{array}{lll}
\add(\HM(\ZZ^\om))& = & \om_1,\\
\cov(\HM(\ZZ^\om))& = & \cov(\iM),\\
\non(\HM(\ZZ^\om))& = & \non(\iM),\\
\cof(\HM(\ZZ^\om))& = &\mathfrak{c}.
\end{array}
\]
\et

\br
In fact, additivity and cofinality works for all non-locally compact Polish groups admitting a two-sided invariant metric, and covering number and uniformity works for non-locally compact Polish groups $G$ admitting a continuous surjective homomorphism onto a non-discrete locally compact Polish group. This latter holds e.g. for $\ZZ^\om$ and for separable infinite dimensional Banach spaces (indeed, $\ZZ^\om$ admits a continuous homomorphism onto $(\ZZ  / 2\ZZ)^\om$, and Banach spaces admit continuous homomorphisms onto their finite dimensional subspaces).
\er

\section{Proofs}

\subsection{Covering number and uniformity}

The results of this section build heavily on \cite{Ba}, in fact, most results of the section are already present in Banakh's paper in some form. However, one of the key points of the present paper is the sharp contrast between the cardinal invariants of $\HN$ and $\HN_{gen}$, so we need to be careful and repeat many familiar argument using this more restrictive definition of Haar nullness.

Although the next lemma is known for the case of the Haar null ideal \cite[Proposition 8]{D}, we include a proof, since our proof for the Haar meager case works for the Haar null case as well.

\bl
Let $\varphi : G \to H$ be a continuous surjective homomorphism between Polish groups. Then the preimage of a Haar null (resp. Haar meager) set is Haar null (resp. Haar meager).
\el

\bp

Applying  \cite[Theorem 4.3 \& Proposition 5.1]{BGJS}  and then \cite[Theorem 11.7]{BGJS} essentially yields the result (here we follow the numbering of the theorems and propositions as in Version 4 of \cite{BGJS}). However, the cited paper only deals with abelian groups, hence, for the convenience of the reader, we include the modified version of \cite[Theorem 11.7]{BGJS} here. (The first two cited results are completely straightforward to adapt to the non-abelian case.)

So let $A \su H$ be Haar null (resp. Haar meager). Without loss of generality we may assume that $A$ is Borel (and hence so is $\varphi^{-1}(A)$). Then by the above cited two results there is a witness function $f : 2^\om \to H$ for $A$. By \cite[Theorem 1.2.6]{BK} $\varphi$ is open (and continuous and onto), hence the multi-function $\Phi : 2^\om \to \iP(G)$ defined as $\Phi(x) = \varphi^{-1} (f(x))$ $(x \in 2^\om)$ is lower semi-continuous (and closed-valued and non-empty-valued).
Therefore by the zero-dimensional Michael Selection Theorem \cite[Theorem 2]{Mi} there is a continuous selection, which in this case means a continuous function $s : 2^\om \to G$ such that $\varphi \circ s = f$. We claim that then $s$ is a witness function for $\varphi^{-1}(A)$. Indeed, an easy calculation shows that $g\varphi^{-1}(A)g' = \varphi^{-1}( \varphi(g) A \varphi(g'))$, and hence
\[
s^{-1} (g\varphi^{-1}(A)g') = s^{-1} ( \varphi^{-1}( \varphi(g) A \varphi(g')) ) = (\varphi \circ s)^{-1} (\varphi(g) A \varphi(g')) = f^{-1 } (\varphi(g) A \varphi(g')),
\]
which is null (resp. meager) in $2^\om$ since $f$ was a witness function for $A$.
\ep

Also recall that in locally compact groups Haar null sets agree with the sets of Haar measure zero and Haar meager sets agree with the meager sets \cite{EN}. Moreover, it is also well-known that $\cov(\iN)$ and $\non(\iN)$ do not depend on the underlying space as long as it is a Polish space equipped with a continuous $\si$-finite Borel measure, e.g. a locally compact non-discrete Polish group equipped with the left Haar measure (where of course we mean that $\iN$ stands for the $\si$-ideal of null sets of the above measure on the Polish space). Similarly, $\cov(\iM)$ and $\non(\iM)$ do not depend on the underlying space as long as it is a Polish space without isolated points, e.g. a non-discrete Polish group (of course  $\iM$ stands for the meager ideal of the Polish space) \cite{BJ}.

\bcor
\lab{c:factor}
If $G$ is a Polish group admitting a continuous surjective homomorphism onto a non-discrete locally compact Polish group then
\[
\begin{array}{lll}
\cov(\HN(G))& \le & \cov(\iN),\\
\cov(\HM(G))& \le & \cov(\iM),\\
\non(\HN(G))& \ge & \non(\iN),\\
\non(\HM(G))& \ge &\non(\iM).
\end{array}
\]
\ecor

\bp
Let $\varphi : G \to H$ be a continuous surjective homomorphism onto a non-discrete locally compact Polish group $H$. By the above remarks, $H$ can be covered by $\cov(\iN)$ many sets of Haar measure zero, which are Haar null in this case, since $H$ is locally compact. Similarly, $H$ can be covered by $\cov(\iM)$ many Haar meager sets.
But then the preimages under $\varphi$ of these sets clearly cover $G$, and these preimages are Haar null (resp. Haar meager) by the previous lemma, finishing the proof of the first two inequalities.

Similarly, the uniformity of the Haar null (resp. Haar meager) sets in $H$ is $\non(\iN)$ (resp. $\non(\iM)$). Choose $X \su G$ with $X \notin \HNG$ and $|X| = \non(\HN(G))$ (resp. $X \notin \HMG$ and $|X| = \non(\HM(G))$). Then $|\varphi(X)| \le |X|$ and $\varphi(X)$ is not of Haar measure zero, or equivalently, $\varphi(X) \notin \HN(H)$ (resp. $\varphi(X) \notin \HM(H)$), otherwise $X \su \varphi^{-1}(\varphi(X)) \in \HNG$ (resp. $\HMG$)  by the previous lemma, a contradiction. 
\ep

Recall that a set $A \su G$ is called \emph{o-bounded} (in symbols $A \in o\iB(G)$) if for each sequence $\{U_n\}_{n \in \om}$ of neighborhoods of the identity there is a sequence $\{F_n\}_{n \in \om}$ of finite subsets of $G$ such that $A \su \cup_{n \in \om} F_n U_n$.

\bl
If $G$ is a non-locally compact Polish group admitting a two-sided invariant metric then $o\iB(G) \su \HNG$.
\el

\bp
This is essentially \cite[Lemma 4]{Ba}, just note that when this paper proves $o\iB(G) \su \HN_{gen}(G)$ the constructed set is in fact Borel.
\ep

\bcor
\lab{c:bd}
If $G$ is a non-locally compact Polish group admitting a two-sided invariant metric then
\[
\begin{array}{lll}
\cov(\HN(G))& \le & \mathfrak{b},\\
\non(\HN(G))& \ge & \mathfrak{d}.
\end{array}
\]
\ecor

\bp
By $o\iB(G) \su \HNG$ we have $\cov(\HN(G)) \le \cov(o\iB(G))$, and by \cite[Lemma 1]{Ba} we have $\cov(o\iB(G)) \le \mathfrak{b}$. Dually, by $o\iB(G) \su \HNG$ we have $\non(\HN(G)) \ge \non(o\iB(G))$, and by \cite[Lemma 1]{Ba} we have $\non(o\iB(G)) \ge \mathfrak{d}$.
\ep

Now we are ready to prove the main results of the section.

\bt

\[
\begin{array}{lll}
\cov(\HN(\ZZ^\om))& = & \min\{\mathfrak{b},  \cov(\iN)\},\\
\non(\HN(\ZZ^\om))& = & \max\{ \mathfrak{d, \non(\iN)}\}.
\end{array}
\]
\et

\bp
Since $\HN(\ZZ^\om) \su \HN_{gen}(\ZZ^\om)$ we obtain $\cov(\HN(\ZZ^\om)) \ge \cov(\HN_{gen}(\ZZ^\om)) = \min\{\mathfrak{b},  \cov(\iN)\}$ and 
$\non(\HN(\ZZ^\om)) \le \non(\HN_{gen}(\ZZ^\om)) = \max\{ \mathfrak{d, \non(\iN)}\}$ using Theorem \ref{t:gen}. The opposite inequalities follow from Corollaries \ref{c:factor} \& \ref{c:bd}.
\ep

\bt
If $G$ is a Polish group admitting a continuous surjective homomorphism onto a non-discrete locally compact Polish group then
\[
\begin{array}{lll}
\cov(\HM(G))& = & \cov(\iM),\\
\non(\HM(G))& = & \non(\iM).
\end{array}
\]
\et

\bp
It is well-known that $\HMG \su \iM$ \cite{EN}, which implies $\cov(\HM(G)) \ge \cov(\iM)$ and $\non(\HM(G)) \le \non(\iM)$. The opposite inequalities follow from Corollary \ref{c:factor}. 
\ep

\subsection{Cofinality}

The main goal of this section is to prove the following.

\bt\lab{t:cof}
If $G$ is a non-locally compact Polish group admitting a two-sided invariant metric then $\cof(\HN(G)) = \cof(\HM(G)) =\mathfrak{c}$.
\et

We will need the following definitions.

\bd
A set $A \su G$ is \emph{compact catcher} if for every compact set $K \su G$ there are $g,h \in G$ such that $gKh \su A$. Let us say that $A$ is \emph{left compact catcher} if for every compact set $K \su G$ there exists $g \in G$ such that $gK \su A$.
\ed

It is easy to see that if $A$ is compact catcher then it is neither Haar null nor Haar meager, see e.g. \cite{EN}. Clearly, the same holds for left compact catcher sets, since left compact catcher sets are compact catcher.

\bigskip
Recall that for a Polish space $G$ the \emph{Effros standard Borel space} of $G$ is denoted by $\iF(G)$. This space consists of the non-empty closed subsets of $G$, and the Borel structure on it is the $\sigma$-algebra generated by the sets of the form $\{F \in \iF(G) : F \cap U \neq \emptyset\}$, where $U \su G$ is open.  

\bigskip
The main technical tool will be the following. 

\bd
We say that the Polish group $G$ is \emph{nice} if there exists a Borel map $\varphi: 2^\om \to \iF(G)$ such that
\begin{enumerate}
\item
if $x \in 2^\om$ then $\varphi(x) \in \HNG \cap \HMG$,
\item
if $P \su 2^\om$ is non-empty perfect then $\cup_{x \in P} \varphi(x)$ is left compact catcher.
\end{enumerate}
\ed

Theorem \ref{t:cof} will immediately follow from the next two results.

\bt\lab{t:nice1}
Every non-locally compact Polish group admitting a two-sided invariant metric is nice.
\et

\bt\lab{t:nice2}
If $G$ is a nice Polish group then $\cof(\HN(G)) = \cof(\HM(G)) =\mathfrak{c}$.
\et

Since the main technical difficulty lies in the proof of Theorem \ref{t:nice1}, we prove Theorem \ref{t:nice2} first. 

\bl\lab{l:coanal}
Let $B \su G$ be a Borel set. Then
\[
\{x \in 2^\om : \varphi(x) \su B \}
\textrm{ is coanalytic}.
\]
\el

\bp
$\{x \in 2^\om : \varphi(x) \su B \} = 
\{x \in 2^\om : \forall y \in G \ (y \in \varphi(x) \implies y \in B ) \} = 
\{x \in 2^\om : \forall y \in G \ (y \notin \varphi(x) \vee y \in B ) \},$
which is a coprojection, hence it suffices to check that
$\{ (x, y) \in 2^\om \times G : y \notin \varphi(x)\}$ and $\{(x, y) \in 2^\om \times G : y \in B \}$ are coanalytic. The latter one is clearly Borel, so we just need to check that the complement of the former one, $\{ (x, y) \in 2^\om \times G : y \in \varphi(x)\}$, is analytic. Indeed,
$\{ (x, y) \in 2^\om \times G : y \in \varphi(x)\} = 
\{ (x, y) \in 2^\om \times G : \exists F \in \iF(G) \ (\varphi(x) = F \wedge y \in F)\}$. Since this is a projection, it suffices to check that the sets $\{ (x, y, F) \in 2^\om \times G \times \iF(G) : \varphi(x) = F \}$ and $\{ (x, y, F) \in 2^\om \times G \times \iF(G) :  y \in F\}$ are analytic. In fact, these sets are already Borel. Indeed, the former one is the graph of a Borel function multiplied by $G$, hence Borel. As for the latter one,
$\{ (x, y, F) \in 2^\om \times G \times \iF(G) :  y \in F\} = 
2^\om \times \{ (y, F) \in  G \times \iF(G) :  y \in F\} = 2^\om \times \cap_n \{ (y, F) \in  G \times \iF(G) : y \in U_n \implies F \cap U_n \neq \emptyset\}$, where $\{U_n\}_{n \in \om}$ is a countable basis of $G$, and this set is Borel since $\{ (y, F) \in  G \times \iF(G) : y \in U_n \}$ is clearly open, and $\{ (y, F) \in  G \times \iF(G) : F \cap U_n \neq \emptyset\} =
G \times \{ F \in \iF(G) : F \cap U_n \neq \emptyset\} $ is Borel by the definition of the Effros space.
\ep

\bigskip
\bp
(Theorem \ref{t:nice2})

First we show that it suffices to prove that for every Haar null (resp. Haar meager) Borel set $B \su G$ we have
\beq\lab{e:cof}
|\{ x \in 2^\om : \varphi(x) \su B \}| \le \om_1.
\eeq
Indeed, first note that if the Continuum Hypothesis holds than we are clearly done, since less than $\mathfrak{c}$ many, in other words countably many nullsets cannot be cofinal, since they cannot even cover $G$, so even a suitable singleton will show that the family is not cofinal. Otherwise, assume that $\ka < \mathfrak{c}$ and $\{ B_\al : \al < \ka \}$ is a cofinal family of Haar null (resp. Haar meager) sets, then by the definition of Haar null (resp. Haar meager) sets without loss of generality we can assume that these sets are Borel. Then every $B_\al$ can contain $\varphi(x)$ for at most $\om_1$ many $x \in 2^\om$, hence there are at most $\ka \cdot \om_1 < \mathfrak{c}$ many $x \in 2^\om$ such that $\varphi(x)$ is contained in some of the $B_\al$'s, contradicting that the family was cofinal.

Now we prove \eqref{e:cof}. Let $B$ be a Haar null (resp. Haar meager) Borel set. Assume to the contrary that $|\{ x \in 2^\om : \varphi(x) \su B \}| > \om_1$. This set is coanalytic by Lemma \ref{l:coanal}, and it is well-known that a coanalytic set of cardinality greater than $\om_1$ contains a non-empty perfect set $P$. (Indeed, this easily follows from the facts that every coanalytic set is the union of $\om_1$ many Borel sets, and every uncountable Borel set contains a non-empty perfect set.) But then $\cup_{x \in P} \varphi(x) \su B$, where $\cup_{x \in P} \varphi(x)$ is non-Haar null (resp. non-Haar meager) by the definition of niceness, and $B$ is Haar null (resp. Haar meager), a contradiction.
\ep

\bigskip
It remains to prove Theorem \ref{t:nice1}.

 \bp
 (Theorem \ref{t:nice1})

 First, it is easy to see that a closed Haar null set is Haar meager. (Indeed, just restrict the witness measure to a compact set such that each relatively open non-empty subset of this set is of positive measure.) Therefore instead of checking $\varphi(x) \in \HNG \cap \HMG$ simply  $\varphi(x) \in \HNG$ will suffice.

For the proof we will need a lot of preparation.

For a finite sequence $s \in S^n \su S^{<\omega}$, $|s|$ will denote $n$, the length of $s$.


The natural numbers are the finite ordinals and we can consider them as von Neumann ordinals, i.e. for every $n \in \omega$
\[ n = \{0,1, \dots, n-1 \}. \]

The following notion will have great importance.
\bdd \label{Tdef}
	A function $\psi$ with $\dom(\psi) = \bigcup_{l \leq i} \prod_{j <l} n_j$ for some $i \in \omega$, $n_0, n_1, \dots, n_{i-1} \in \omega \setminus \{0 \}$ is an element of $\mathcal{T}$ (and we call $\psi$ a ''labeled tree''), if
	$\psi: \dom(\psi) \to 2^{<\omega}$ is an injective mapping, such that
	\begin{equation} \label{f1} \psi(\emptyset) = \emptyset, \end{equation}
	and
	\begin{equation} \label{f2} s \su t \ \iff \ \psi(s) \su \psi(t), \end{equation}
	i.e. $\psi$ maps end-extensions to end-extensions, and incomparable elements to incomparable elements.
	(Under $\prod_{j<0} n_j$ we mean the set containing exactly the empty sequence $\emptyset $, thus $\emptyset \in \dom(\psi)$. 
	Moreover, for the function $\pi$ defined only on $\emptyset$, mapping it to $\emptyset$, we have $\pi \in \mathcal{T}$.)
\ed

\bdd \label{nDef}
	For  $\psi \in \mathcal{T}$, define $\underline{n}^\psi \in (\omega \setminus \{0 \})^{<\omega}$ so that
	\[ \dom(\psi) = \bigcup_{l \leq i} \prod_{j <l} n^\psi_j. \]
	
\ed


Note that $\mathcal{T}$ is countable.

We partition $\mathcal{T}$ as follows.
\bdd\label{T__k} 
	\begin{equation} \label{T_k} \mathcal{T}_k = \{ \psi \in \mathcal{T}: \ \underline{n}^\psi \in \omega^k \}. \end{equation}
	
\ed
Then
\[ \mathcal{T} = \bigcup_{k \in \omega} \mathcal{T}_k. \]

\br \label{mcT_0}
	Since $\psi \in \mathcal{T}$ maps $\emptyset$ to $\emptyset \in 2^{<\omega}$ (Definition $\ref{Tdef}$)
	\[  \{ \langle \emptyset,  \emptyset \rangle \} \in \mathcal{T}_0,\]
	that is the unique element of $\mathcal{T}_0$.
	
\er

Now we define a partial order on $\mathcal{T}$.

\bdd \label{Tpo}
	$\mathcal{T}$ is a poset with the following partial order:
	\[ \psi \leq \psi' \]
	iff  
	\begin{enumerate}
		\item \label{POf1} for the sequences $\underline{n}^\psi, \underline{n}^{\psi'} \in \omega^{<\omega}$,  $\underline{n}^{\psi'}$ is an end-extension of $\underline{n}^\psi$, i.e.
		\[\underline{n}^\psi \su \underline{n}^{\psi'}, \]
		
		\item \label{POf2} $\psi \su  \psi'$.
	\end{enumerate}
\ed

	In particular if $\psi \leq \psi'$, then  $\psi = (\psi')\restrict{\dom(\psi)}$ and for the length $l= |\underline{n}^{\psi}|$ of $\underline{n}^{\psi} \in \omega^{<\omega}$ we have $\dom(\psi) = \dom(\psi') \cap \omega^{\leq l}$.

	Combining this with Remark $\ref{mcT_0}$ we have the following.
	\br \label{megsz0}
		$\psi \geq \{ \langle \emptyset , \emptyset \rangle \}$ for every $\psi \in \mathcal{T}$.
	\er

	\bfa \label{Tfa}
		Suppose that $\psi \in \mathcal{T}_i$. Then for $k \leq i$ defining 
		\begin{equation} \label{psik} \psi_k = \psi\restrict{\dom(\psi) \cap \omega^{\leq k}} \end{equation} 
		we have
		$\psi_k \in \mathcal{T}_k$, moreover,
		$\psi_0 = \{ \langle \emptyset, \emptyset \rangle \} \leq \psi_1 \leq \dots \leq \psi_{i-1} \leq \psi_i = \psi$
		are the only elements less than or equal to $\psi$ in $\mathcal{T}$ (in the sense of Definition $\ref{Tpo}$).
		
		In particular, if $\psi \neq \pi \in \mathcal{T}_i$, then $\psi$ and $\pi$ are incomparable w.r.t. the partial order on $\mathcal{T}$.
	\efa
	

	\bdd \label{mdef}
		We define an embedding of $(\mathcal{T}, \leq )$ into $(\omega^{<\omega}, \su )$ as follows. To each
		 $\psi \in \mathcal{T}$ we assign by induction on $|\underline{n}^\psi|$ a sequence $\underline{m}^\psi \in \omega^{<\omega}$ such that
		\begin{equation}
		|\underline{m}^\psi| = |\underline{n}^\psi|, \text{ i.e. } \underline{m}^\psi \in \omega^{|\underline{n}^\psi|}.
		\end{equation}
	
	\ed

 First, fix a compatible two-sided invariant metric $d$, i.e. for which
 \begin{equation} \label{invd} d(x,y) = d(zx,zy)= d(xz,yz) \ \forall x,y,z \in G. \end{equation} 
 (Such a metric is also automatically complete \cite[Corollary 1.2.2]{BK}.)

 \bdd
 	For $r>0$ let $B(r) \su G$ denote the ball of radius $r$ centered at the identity of the group $e$, i.e.
 	\[ B(r) = \left\{ g \in G \ : \ d(g,e) < r \right\}. \]
 \ed
 
 The following two lemmas state well-known facts using the invariance of $d$, we leave the proof to the reader.
 \bl \label{szorzasfact}
 	Let $g,h \in G$. Then $d(gh,e) \leq d(g,e)+d(h,e)$, moreover, for each $\varepsilon,\varepsilon'  >0$ we have $B(\varepsilon)B(\varepsilon') \su B(\varepsilon+ \varepsilon')$.
 	 \el

 \bl \label{gomb}
 	If $g \in G$, $\varepsilon>0$, then for the open ball $B(\varepsilon)$
 	\[ g B(\varepsilon) = B(\varepsilon) g = \{h \in G: \ d(g,h) < \varepsilon\}. \]
 \el

 The next technical step is similar to the one in the proof of the main theorem in \cite{solecki} used for constructing $2^\omega$-many pairwise disjoint compact catcher sets.
 \bl \label{fo}
 	Let $G$ be a Polish non-locally compact group with the two-sided invariant metric $d$. Then there exist sequences $(\delta_i)_{i \in \omega}$, 
 	$(\varepsilon_i)_{i \in \omega}$, $\left(q^{(i)}_j\right)_{j \in \omega}$ ($i \in \omega$) such that
 	\begin{enumerate}[(i)]
 		\item $q^{(i)}_0 =e$ for each $i \in \omega$,
 		\item \label{gk} if $i >0$ then $\{q^{(i)}_j \ : \ j \in \omega \} \su B\left(\frac{\delta_{i-1}}{2}\right)$, 
 		\item \label{suru2}
 		\begin{enumerate}[(a)]
 			\item  $\cl(\{q^{(i)}_j \ : \ j \in \omega \}) = \cl(B(\frac{\delta_{i-1}}{2}))$ if $i>0$, i.e. it is a countable dense subset of the ball,
 			\item $\cl(\{q^{(0)}_j \ : \ j \in \omega \}) = G$,
 		\end{enumerate}
 		\item \label{del1} $\delta_{i} + \varepsilon_{i} < \frac{\delta_{i-1}}{2}$ for each $i >0$,
 		\item  \label{del15} $d(q^{(i)}_0, q^{(i)}_1) = 5 \delta_i$,
 		\item \label{del2} $B(\varepsilon_i)$ cannot be covered by finitely many subsets of diameter at most $20 \delta_i$.
 		
 	\end{enumerate}
 \el
 \bp
 	Let $\varepsilon_0 = 1$, and choose a countable dense subset $\{q^{(0)}_0 = e, q^{(0)}_1, \dots, q^{(0)}_j, \dots \} \su G$ so that there is no finite cover of $B(\varepsilon_0) = B(1)$  with sets of diameter at most $4\cdot d(q^{(0)}_{0}, q^{(0)}_{1})$, and define $\delta_0 = \frac{d(q^{(0)}_{0}, q^{(0)}_{1})}{5}$.
 	
 	Fix $i>0$, and assume that $\delta_k$, $(q^{(k)}_j)_{j \in \omega}$, $\varepsilon_k$ are already defined for each $k<i$.
 	Let $\varepsilon_i = \frac{\delta_{i-1}}{4}$. There exists a constant $c>0$ such that there is no finite $c$-net in $B(\varepsilon_i)$ (hence $c < \varepsilon_i$).
 	Now if $g \in B(\frac{\delta_{i-1}}{2})$ is such that $d(g,e) \leq \frac{c}{4} $, then choosing $q^{(i)}_1=g$, $q^{(i)}_0 =e$, $\delta_i = d(g,e)/5 \leq \frac{c}{20} < \frac{\varepsilon_i}{20} = \frac{\delta_{i-1}/4}{20}$
 	will satisfy conditions $\eqref{del1}$, $\eqref{del2}$. Also,
 	the set $\{ q^{(i)}_2,q^{(i)}_3, \dots, q^{(i)}_n, \dots \} \su B(\frac{\delta_{i-1}}{2})$ can be chosen so that  condition $\eqref{suru2}$ holds.
 \ep

 \bdd \label{Qdef}
 	For each $i \in \omega$ and $j>0$ let 
 	\[ Q^{(i)}_j= \{q^{(i)}_0,q^{(i)}_1, \dots ,q^{(i)}_{j-1} \}. \]
 \ed
 By our previous lemma we can deduce the following.
 \bl \label{eltolos}
 	For each $i \in \omega$ there are elements $t^{(i)}_{n,m} \in B(\varepsilon_i)$ ($n,m \in \omega$, $n>0$) such that
 	whenever $\langle n,m \rangle \neq \langle n',m' \rangle$
 	\begin{equation} \label{tav} d( Q^{(i)}_{n}t^{(i)}_{n,m}, Q^{(i)}_{n'} t^{(i)}_{n',m'}) > 9 \delta_i. \end{equation}
 \el
 \bp
 	Fix $i \in \omega$, and fix an enumeration of $\omega \times \omega$:
 	\[ \omega \times \omega = \{ \langle n_k,m_k \rangle : \ k \in \omega \}. \]
 	Let $j \in \omega$ $j>0$ be fixed, and assume that the $t^{(i)}_{n_l,m_l}$'s are determined for $l <j$. Let $Q = Q^{(i)}_{n_j} = \{q^{(i)}_0, q^{(i)}_1, \dots, q^{(i)}_{n_j-1}\}$.
 	
 	Now assume on the contrary that there is no $t \in B(\varepsilon_i)$ such that
 	\[ \forall l<j: \ d( Qt, Q^{(i)}_{n_l}t^{(i)}_{n_l,m_l} ) > 9 \delta_i. \]
 	This means that for each $t \in B(\varepsilon_i)$ there is an element $q \in Q$, and $l<j$ such that
 	\begin{equation} \label{suruen}  d( qt,  Q^{(i)}_{n_l} t^{(i)}_{n_l,m_l}) \leq 9 \delta_i,\end{equation}
 	thus by the invariance of $d$
 	\[ d(t, q^{-1} Q^{(i)}_{n_l} t^{(i)}_{n_l,m_l}) \leq   9 \delta_i. \]
 	Using Lemma $\ref{gomb}$
 	\[ t \in q^{-1}  Q^{(i)}_{n_l} t^{(i)}_{n_l,m_l} B(10 \delta_i) = \bigcup_{q' \in  Q^{(i)}_{n_l}} q^{-1} q' t^{(i)}_{n_l,m_l}  B(10 \delta_i). \]
 	And because for each $t \in B(\varepsilon_i)$ there exist   $q \in Q$, and $l<j$ such that $\eqref{suruen}$ holds, we obtain that
 	\[ B(\varepsilon_i) \su \bigcup_{q \in Q, l < j} q^{-1}  Q^{(i)}_{n_l} t^{(i)}_{n_l,m_l} B(10 \delta_i) = \]
 	\[ = \bigcup_{q \in Q, l < j} \bigcup_{q' \in Q^{(i)}_{n_l}} q^{-1} q'  t^{(i)}_{n_l,m_l} B(10 \delta_i). \]
 	Using the finiteness of $Q = Q^{(i)}_{n_j}$ and the $Q^{(i)}_{n}$'s (Definition $\ref{Qdef}$) and recalling Lemma $\ref{gomb}$ we conclude that $B(\varepsilon_i)$ can be covered by finitely many subsets of diameter at most $20 \delta_i$, contradicting $\eqref{del2}$ from Lemma $\ref{fo}$.
 \ep


\bdd \label{gdef}
	For every integer $i>0$ and sequences $p \in (\omega \setminus \{0\})^i$, $r, s \in \omega^{i}$ we define the group element $g_{p,r,s} \in G$
	\[ g_{p,r,s} = \left( q^{(0)}_{s_{0}} t^{(0)}_{p_{0},r_{0}} \right)  \left(q^{(1)}_{s_{1}} t^{(1)}_{p_{1},r_{1}} \right) \ldots  \left(  q^{(i-1)}_{s_{i-1}} t^{(i-1)}_{p_{i-1},r_{i-1}} \right). \]
\ed

Lemma $\ref{eltolos}$ yields the following. 
\bl \label{tavlemma}
	Suppose that $i>0$, $\underline{n} = (n_j)_{j<i},$ $\underline{n}' = (n'_j)_{j<i}  \in (\omega \setminus \{0\})^i$ and $\underline{m} = (m_j)_{j<i},$ $\underline{m}' = (m'_j)_{j <i} \in \omega^i$ are such that
	\[ (\langle n_j,m_j \rangle)_{j <i} \neq (\langle n_j',m_j' \rangle)_{j < i}, \]
	and they differ only on the last, $i-1$-th coordinate, i.e.
	\[ (\langle n_j,m_j \rangle)_{j <i-1} = (\langle n_j',m_j' \rangle)_{j < i-1}, \]
	 
	Let $s \in \prod_{j<i} n_j$, $s' \in \prod_{j<i} n_j'$ be such that $s\restrict{i-1} = (s')\restrict{i-1}$. Then
	\begin{equation}
	d(g_{\underline{n},\underline{m},s},g_{\underline{n}',\underline{m}',s'}) > 9 \delta_{i-1}.
	\end{equation}
\el
\bp
	Since the two products have the same first $i-1$ coordinates, i.e.
	\[
	(q^{(0)}_{s_{0}} t^{(0)}_{n_{0},m_{0}} )(q^{(1)}_{s_{1}} t^{(1)}_{n_{1},m_{1}} ) \dots ( q^{(i-2)}_{s_{i-2}}t^{(i-2)}_{n_{i-2},m_{i-2}})=  \]
	\[ ( q^{(0)}_{s'_{0}}t^{(0)}_{n'_{0},m'_{0}})( q^{(1)}_{s'_{1}}t^{(1)}_{n'_{1},m'_{1}}) \dots ( q^{(i-2)}_{s'_{i-2}}t^{(i-2)}_{n'_{i-2},m'_{i-2}}),\]
	by the invariance of the metric $d$ it is enough to prove that  
	\begin{equation} \label{cld}
	d\left(q^{(i-1)}_{s_{i-1}}t^{(i-1)}_{n_{i-1},m_{i-1}}, q^{(i-1)}_{s'_{i-1}}t^{(i-1)}_{n'_{i-1},m'_{i-1}}\right) > 9 \delta_{i-1}.
	\end{equation}
	Now, since $s_{i-1} \in n_{i-1}$, and $s'_{i-1} \in n'_{i-1}$ by our assumptions, using Definition $\ref{Qdef}$ $q^{(i-1)}_{s_{i-1}} \in Q^{(i-1)}_{n_{i-1}}$ and  $q^{(i-1)}_{s'_{i-1}} \in Q^{(i-1)}_{n'_{i-1}}$. This yields using Lemma $\ref{eltolos}$ and $\langle n_{i-1},m_{i-1} \rangle) \neq \langle n'_{i-1},m'_{i-1} \rangle$ that
	\begin{equation} \label{fod}
	d\left(  q^{(i-1)}_{s_{i-1}}t^{(i-1)}_{n_{i-1},m_{i-1}},  q^{(i-1)}_{s'_{i-1}} t^{(i-1)}_{n'_{i-1},m'_{i-1}} \right) > 9 \delta_{i-1}.
	\end{equation}
\ep

Now we can turn to the construction of $\varphi$ from Proposition $\ref{t:nice1}$. Fix $c \in 2^\omega$.
\bdd \label{cfaj}
	For $c \in 2^\omega$ let 
	\begin{equation} \label{cfaja} \mathcal{T}(c) = \{ \psi \in \mathcal{T} : \ \exists t \in \prod_{j<|\underline{n}^\psi|} n^\psi_j: \ \psi(t) \su c \}, \end{equation}
	i.e. there is an element
	\[ t \in \prod_{j<|\underline{n}^\psi|} n^\psi_j \su \dom(\psi) = \bigcup_{l \leq |\underline{n}^\psi|} \prod_{j<l} n^\psi_j \]
	for which $c$ is an end-extension of $\psi(t)$.
	
	For $k \in \omega$ define
	\begin{equation} \label{ckfaja}
	\mathcal{T}_k(c) = \mathcal{T}(c) \cap \mathcal{T}_k. \end{equation} 
\ed
Note that condition $\eqref{f2}$ from Definition $\ref{Tdef}$ implies that (for a fixed $\psi \in \mathcal{T}$) at most one 
$t \in \prod_{j<|\underline{n}^\psi|} n^\psi_j$ can exist for which $\psi(t) \in 2^{<\omega}$ is an initial segment of $c$.

\bdd \label{ccdef}
	For $c \in 2^\omega$, $\psi \in \mathcal{T}(c)$  let $b^{c,\psi} \in \prod_{j<|\underline{n}^\psi|} n^\psi_j \su \dom(\psi)$ be the element for which
	\[ \psi(b^{c,\psi}) \su c. \]
\ed

\br\label{mcT_0(c)}
	\[ \{ \langle \emptyset, \emptyset  \rangle \} \in \mathcal{T}_0(c) = \mathcal{T}_0  \ \ (\forall \ c \in 2^\omega). \]
\er

The following is an easy observation, the proof is left to the reader.
\bl \label{cmegsz}
	Let $c \in 2^\omega$, let $\psi \in \mathcal{T}_i(c) \su \mathcal{T}(c)$. Then for any $\pi \leq \psi$ (in the sense of Definition $\ref{Tpo}$) with  $\pi \in \mathcal{T}_k$ we have
	\begin{equation} \label{T(c)} \pi  \in \mathcal{T}_{k}(c),\end{equation}
	moreover,
	\begin{equation} \label{cbeta} b^{c,\pi} = (b^{c,\psi})\restrict{k}. \end{equation}
\el



\bdd \label{Fdef}
	Let $\psi \in \mathcal{T}_i(c)$, then (recall the definitions of $\underline{n}^\psi$, $\underline{m}^\psi$, Definitions $\ref{nDef}$, $\ref{mdef}$,)  for $i>0$ define
	\begin{equation}
	\begin{array}{l}
	F_c^\psi =\cl(B(\delta_{i-1})) g_{\underline{n}^\psi,\underline{m}^\psi,b^{c,\psi}} = \\
	= \cl(B(\delta_{i-1}))  \left( q^{(0)}_{b^{c,\psi}_{0}} \cdot t^{(0)}_{n^\psi_{0},m^\psi_{0}} \right)  \left(q^{(1)}_{b^{c,\psi}_{1}} \cdot t^{(1)}_{n^\psi_{1},m^{\psi}_{1}} \right) \ldots  \left(  q^{(i-1)}_{b^{c,\psi}_{i-1}} \cdot t^{(i-1)}_{n^\psi_{i-1},m^{\psi}_{i-1}} \right).
	\end{array}
	\end{equation}
	
	For $i=0$ and the unique element $\psi = \{  \langle \emptyset, \emptyset \rangle \} $
	of $\mathcal{T}_0 = \mathcal{T}_0(c)$, define
	\begin{equation} \label{nulla}
	F_c^\psi = G.
	\end{equation}

\ed

We are ready to define $\varphi(c)$.
\bdd \label{varphidef}
	\begin{equation} \label{fc} \begin{array}{l} \varphi(c) = \bigcap_{i \in \omega} \bigcup_{\psi \in \mathcal{T}_i(c)} F^\psi_c = \bigcap_{i >0} \bigcup_{\psi \in \mathcal{T}_i(c)} \cl(B(\delta_{i-1})) g_{\underline{n}^\psi,\underline{m}^\psi,b^{c,\psi}} =\\
	= \bigcap_{i > 0} \bigcup_{\psi \in \mathcal{T}_i(c)} \cl(B(\delta_{i-1}))  \left( q^{(0)}_{b^{c,\psi}_{0}} \cdot t^{(0)}_{n^\psi_{0},m^\psi_{0}} \right)  \left(q^{(1)}_{b^{c,\psi}_{1}} \cdot t^{(1)}_{n^\psi_{1},m^{\psi}_{1}} \right) \ldots  \left(  q^{(i-1)}_{b^{c,\psi}_{i-1}} \cdot t^{(i-1)}_{n^\psi_{i-1},m^{\psi}_{i-1}} \right). \end{array} \end{equation}
\ed
The following lemmas will ensure that $\varphi(c)$ is closed.

\bl \label{monotonGomb}
	Let $i \in \omega$, $p \in (\omega \setminus \{0\})^i$, $r \in \omega^i$, $s \in \prod_{j<i} p_j$ and $k\leq i$ be fixed.
	Then
	\[ \cl(B(\delta_{k-1}))g_{p\restrict{k}, r\restrict{k}, s\restrict{k}} \supset \cl(B(\delta_{i-1})) g_{p,r,s}.\]
\el
\bp
	W.l.o.g. we can assume that $k = i-1$. Then by Definition $\ref{gdef}$
	\[ g_{p,r,s} = g_{p\restrict{i-1}, r\restrict{i-1}, s\restrict{i-1}}  q^{(i-1)}_{s_{i-1}} t^{(i-1)}_{p_{i-1},r_{i-1}} . \]
	Hence by the invariance of $d$ it is enough to show that 
	\begin{equation} \label{kozel}  q^{(i-1)}_{s_{i-1}}t^{(i-1)}_{p_{i-1},r_{i-1}} \in B(\delta_{i-2}- \delta_{i-1}), \end{equation} 
	that is, 
	\begin{equation} \label{kozelm} d( q^{(i-1)}_{s_{i-1}}t^{(i-1)}_{p_{i-1},r_{i-1}},e) < \delta_{i-2}- \delta_{i-1}. \end{equation} 
	
	Now $t^{(i-1)}_{p_{i-1},r_{i-1}} \in B(\varepsilon_{i-1})$ by the definition of $t^{(j)}_{n,m}$'s (Lemma $\ref{eltolos}$),
	also  $q^{(i-1)}_{s_{i-1}} \in B(\delta_{i-2}/2)$ by $\eqref{gk}$ in Lemma $\ref{fo}$. Therefore $d(q^{(i-1)}_{s_{i-1}}t^{(i-1)}_{p_{i-1},r_{i-1}}, e)< \varepsilon_{i-1} + \delta_{i-2}/2$ by Lemma $\ref{szorzasfact}$.
	Finally using the fact that $\varepsilon_{i-1} + \delta_{i-1} < \delta_{i-2}/2$ ($\eqref{del1}$ in Lemma $\ref{fo}$)
	\[ d(q^{(i-1)}_{s_{i-1}}t^{(i-1)}_{p_{i-1},r_{i-1}}, e) + \delta_{i-1} < \varepsilon_{i-1} +  \delta_{i-2}/2 + \delta_{i-1} < \delta_{i-2} \]  
	yields $\eqref{kozelm}$.
\ep

\bl \label{monotonF}
	Let $c \in 2^\omega$, $\psi \in \mathcal{T}_i(c)$, $\pi \in \mathcal{T}_k(c)$, assume that $k>0$. Then $\psi \geq \pi$ implies 
	\[ 	F_c^\psi \su F_c^\pi. \]
\el
\bp
	
	Using Corollary $\ref{Tfa}$ w.l.o.g. we can assume that $\pi = \psi\restrict{(\dom(\psi) \cap \omega^{\leq i-1})} \in \mathcal{T}_{i-1}$, i.e. $k = i-1 >0$.
	Since $\pi \leq \psi$, $(\underline{n}^\psi)\restrict{i-1} = \underline{n}^\pi$ (by Definition  $\ref{Tpo}$ $\eqref{POf1}$), and $(\underline{m}^\psi)\restrict{i-1} = (\underline{m}^\pi)$ (by Definition $\ref{mdef}$).
	Claim $\ref{cmegsz}$ gives that $\pi \in \mathcal{T}_k(c)$, and $b^{c,\psi}\restrict{i-1} = b^{c,\pi}$. Finally, applying Lemma $\ref{monotonGomb}$ with $p = \underline{n}^\psi$,
	$r = \underline{m}^\psi$, $s = b^{c,\psi} \in \prod_{j<i} n^\psi_j$ finishes the proof.
\ep

\bl \label{tavF}
	Let $c \in 2^\omega$, $\psi, \pi \in \mathcal{T}_i(c)$ be fixed. Suppose that for $\psi_{i-1} = \psi\restrict{(\dom(\psi) \cap \omega^{\leq i-1})}$
	and $\pi_{i-1} = \pi\restrict{(\dom(\pi) \cap \omega^{\leq i-1})}$ $\psi_{i-1} = \pi_{i-1}$, but $\psi \neq \pi$.  Then the following inequality holds
	\begin{equation} \label{ie} d(g_{\underline{n}^\psi,\underline{m}^\psi,b^{c,\psi}}, g_{\underline{n}^\pi,\underline{m}^\pi,b^{c,\pi}}) > 9 \delta_{i-1}.
	\end{equation}
\el
\bp
	First, applying Lemma $\ref{cmegsz}$ we obtain $\psi_{i-1}= \pi_{i-1} \in \mathcal{T}_{i-1}(c)$ (we used $\psi_{i-1} \leq \psi$, $\pi_{i-1} \leq \pi$ by Claim $\ref{Tfa}$), moreover if $s = b^{c,\psi}$, $s' = b^{c,\pi}$ then 
	\begin{equation} \label{segyenl} s\restrict{i-1} =(b^{c,\psi})\restrict{i-1} = b^{c,\psi_{i-1}} = b^{c,\pi_{i-1}} = (b^{c,\pi})\restrict{i-1} = (s')\restrict{i-1}. \end{equation} 
	By the definition of the partial order on $\mathcal{T}$ (Definition $\ref{Tpo}$), and Definitions $\ref{nDef}$, $\ref{mdef}$ we have that 
	$\left( \langle n^{\psi_{i-1}}_j, m^{\psi_{i-1}}_j \rangle \right)_{j<i-1} = \left( \langle n^\psi_j, m^\psi_j \rangle \right)_{j<i-1}$ and
	$\left( \langle n^{\pi_{i-1}}_j, m^{\pi_{i-1}}_j \rangle \right)_{j<i-1} = \left( \langle n^\pi_j, m^\pi_j \rangle \right)_{j<i-1}$. Therefore, since $\psi_{i-1} = \pi_{i-1}$

	\[ \left( \langle n^\psi_j, m^\psi_j \rangle \right)_{j<i-1} = \left( \langle n^\pi_j, m^\pi_j \rangle \right)_{j<i-1}, \]
		and
		\[ m^\psi_{i-1} \neq m^\pi_{i-1} \]
   		(since $\underline{m}^\psi \neq \underline{m}^\pi$ by $\psi \neq \pi$).
   		Now using $\eqref{segyenl}$, and the fact that $s = b^{c,\psi} \in \prod_{j<i} n^\psi_j$, $s' = b^{c,\pi} \in \prod_{j<i} n^\pi_j$ (by Definition $\ref{ccdef}$), we can apply Lemma $\ref{tavlemma}$ with 
   		$( \langle n_j,m_j \rangle )_{j <i} = (\langle n^\psi_j, m^\psi_j \rangle )_{j<i}$, and
   		$( \langle n'_j,m'_j \rangle )_{j <i} = (\langle n^\pi_j, m^\pi_j \rangle )_{j<i}$. This yields $\eqref{ie}$, as desired.
\ep


\bcor \label{Ffa}
	Let $c \in 2^\omega$, $\psi \in \mathcal{T}_i(c), \pi \in \mathcal{T}_k(c)$.
	Then  $\pi \leq \psi$ implies   $F^\psi_c \su F^\pi_c$,
	and if $\psi$ and  $\pi$ are incomparable in the partial order of $\mathcal{T}$, then 
	\begin{equation} \label{Fektav} d(F^\psi_c,F^\pi_c) >  7 \delta_{\min\{i,k\}-1 }. \end{equation}
\ecor
\bp
	First, we assume that $\pi \leq \psi$ (hence $k \leq i$). If $k>0$, then apply Lemma $\ref{monotonF}$, and we are done.
	If $k = 0$, then $\pi$ is the unique element of $\mathcal{T}_0 = \mathcal{T}_0(c)$ (Remark $\ref{mcT_0(c)}$), and $F^\pi = G$ by Definition $\ref{Fdef}$, thus  $F^\pi_c \supset F^\pi_c$ obviously holds.  
	
	Now assume that $\psi, \pi$ are incomparable in the partial order of $\mathcal{T}$.
	Let $i_0$ be such that $\psi_{i_0-1} = \psi\restrict{(\dom(\psi) \cap \omega^{\leq i_0-1})}$ and $\pi_{i_0-1} = \pi\restrict{(\dom(\pi) \cap \omega^{\leq i_0-1})}$ 
	are equal, but $\psi_{i_0}= \psi\restrict{(\dom(\psi) \cap \omega^{\leq i_0})} \neq \pi_{i_0} = \pi\restrict{(\dom(\pi) \cap \omega^{\leq i_0})}$ (such an $i_0$ exists by Claim $\ref{Tfa}$, also $i_0 \leq i,k$, $\pi_{i_0} \leq \pi$, $\psi_{i_0} \leq \psi$).
	Now $\psi_{i_0}, \pi_{i_0} \in \mathcal{T}_{i_0}(c)$ by Lemma $\ref{cmegsz}$. 
	Applying Claim $\ref{tavF}$ for	$\eta = \psi_{i_0}, \zeta = \pi_{i_0} \in \mathcal{T}_{i_0}(c)$ gives that 
	\[  d(g_{\underline{n}^\eta,\underline{m}^\eta,b^{c,\eta}},g_{\underline{n}^\zeta,\underline{m}^\zeta,b^{c,\zeta}}) > 9 \delta_{i_0-1}.  \]
	This with the definition $F_c^\eta, F_c^\zeta$ (Definition $\ref{Fdef}$) gives that
	\[ d(F_c^\eta, F_c^\zeta)> 7 \delta_{i_0-1}.\]
	Moreover, $\eta \leq \psi$, $\zeta \leq \pi$ implies using Lemma $\ref{monotonF}$ that $F_c^\eta \supset F_c^\psi$, $F_c^\zeta \supset F_c^\pi$, therefore $d(F^\psi_c,F^\pi_c) >  7 \delta_{i_0-1}$.
	Finally, the inequality $i_0 \leq i,k$, and the monotonicity of the $\delta$'s ($\eqref{del1}$ from Lemma $\ref{fo}$) implies $\eqref{Fektav}$.
\ep

It is worth mentioning the following consequence of Lemma $\ref{monotonF}$.
\bcor \label{Ffa2}
	For a fixed $c \in 2^\omega$, 	$i \leq j \in \omega$ implies that  $\bigcup_{\psi \in \mathcal{T}_i(c)} F^\psi_c \supset \bigcup_{\psi \in \mathcal{T}_j(c)} F^\psi_c$.
\ecor

\bl
	Let $c \in 2^\omega$ be given. Then 
	\[ \varphi(c) = \bigcap_{i \in \omega} \bigcup_{\psi \in \mathcal{T}_i(c)} F^\psi_c \]
	is closed.
\el
\bp
	It is enough to show that for a fixed $i \in \omega$ the set
	\[ \bigcup_{\psi \in \mathcal{T}_i(c)} F^\psi_c = \bigcup_{\psi \in \mathcal{T}_i(c)}  \cl(B(\delta_{i-1})) g_{\underline{n}^\psi,\underline{m}^\psi, b^{c,\psi}} \]
	is closed. Recall the fact that for a system of closed sets $\{ F_i: \ i \in \mathcal{I} \}$ with a constant $d_0>0$ such that $d_0 \leq d(F_i, F_j)$
	whenever $i \neq j \in \mathcal{I}$, the union $\bigcup_{i \in \mathcal{I}} F_i$ is closed.
	Therefore (by the invariance of $d$) it suffices to show that there is a constant $r_i> 0$ such that for $\psi \neq \pi \in \mathcal{T}_i(c)$
	\begin{equation} \label{ri}
	r_i \leq d(F^\psi_c, F^\pi_c).
	\end{equation}
	But for fixed $\psi \neq \pi \in \mathcal{T}_i(c)$ applying Corollary $\ref{Ffa}$ we have
	\[ 7 \delta_{i-1} \leq d(F^\psi_c, F^\pi_c),\]
	i.e. we proved $\eqref{ri}$ with $r_i = 7 \delta_{i-1}$.
\ep

\bl
	Let $c \in 2^\omega$. Then $F = \varphi(c)$ is Haar null.
\el
\bp
	First we define a measure $\mu$ according to which each two-sided translate of $\varphi(c)$ will be null. For a fixed $i \in \omega$ and $s \in 2^i$ define
	\begin{equation} \label{Cs} C_s = \left\{ \begin{array}{ll}   \cl(B(\delta_{i-1})) q^{(0)}_{s_0} \dots  q^{(i-2)}_{s_{i-2}}q^{(i-1)}_{s_{i-1}}& \text{ if } i>0, \\
	G & \text{ if } i= 0, \text{ i.e. } s = \emptyset. \end{array} \right.
	\end{equation}
	Let
	\begin{equation} \label{C^i} C^i = \bigcup_{s \in 2^i} C_s, \end{equation}
	and define
	\[ C = \bigcap_{i \in \omega} C^i. \]
	We will prove that $C$ is a Cantor set, i.e. it is homeomorphic to $2^\omega$ by showing that 
	\begin{equation} \label{diszj} \text{if } s \neq s' \in 2^i \ \text{  then }C_s \cap C_{s'} = \emptyset, \end{equation}
	and 
	\begin{equation} \label{csokk} \text{ if }s,t \in 2^{<\omega}, \ s \subset t \ \text{ then }C_s \supset C_t \end{equation}
	(which together with $\eqref{C^i}$ will imply that
	\begin{equation} \label{Cimon}
	i\leq j \ \ \to \ C^i \supset C^j). 
	\end{equation} 
	(Recall that since $\diam(C_s) \leq 2 \delta_{|s|-1}$ by the invariance of $d$, and the $\delta_i$'s converge to $0$ by $\eqref{del1}$ from Lemma $\ref{fo}$, therefore the mapping $2^\omega \ni c \mapsto \cap_{n \in \omega} C_{c\restrict{n}}$ is indeed a homeomorphism, see \cite[(6.2)]{Ke}.)
	
	First, if $i>0$ then since $q^{(i)}_{s_i}$ equals either $q^{(i)}_0$, or  $q^{(i)}_1$, and $q^{(i)}_0,q^{(i)}_1 \in B(\frac{\delta_{i-1}}{2})$ (by $\eqref{gk}$ in Lemma $\ref{fo}$),  thus the invariance of $d$ implies (fixing $s \in 2^{i+1}$)
	\[  d(q^{(0)}_{s_{0}}q^{(1)}_{s_{1}} \dots q^{(i-1)}_{s_{i-1}} \cdot e, q^{(0)}_{s_{0}} q^{(1)}_{s_{1}} \ldots q^{(i-1)}_{s_{i-1}} \cdot q^{(i)}_{s_i}) < \delta_{i-1}/2. \]
	Moreover, the fact that $\delta_{i}< \delta_{i-1}/2$ (by $\eqref{del1}$ from Lemma $\ref{fo}$) clearly implies that (using Lemma $\ref{gomb}$)
	\[B(\delta_{i-1}) q^{(0)}_{s_{0}}q^{(1)}_{s_{1}} \dots q^{(i-1)}_{s_{i-1}}  \supset B(\delta_{i}) q^{(0)}_{s_0}q^{(1)}_{s_{1}} \dots q^{(i-1)}_{s_{i-1}} q^{(i)}_{s_{i}} . \]
	After taking closures we obtain $\eqref{csokk}$.
	
	Now as we already have $\eqref{csokk}$, for $\eqref{diszj}$ it suffices to show that for a fixed $0<i \in \omega$ if $s \neq s' \in 2^i$, and $s\restrict{i-1} = s'\restrict{i-1}$ (i.e. the $i-1$-th is the first coordinate on which $s$ and $s'$ differ) then
	\begin{equation} \label{kell} d\left(  q^{(0)}_{s_0}q^{(0)}_{s_{1}} \ldots q^{(i-1)}_{s_{i-1}}, q^{(0)}_{s'_{0}}q^{(1)}_{s'_{1}} \dots q^{(i-1)}_{s'_{i-1}} \right)> 2 \delta_{i-1}. \end{equation}
	But in this case, since $s_{i-1} \neq s'_{i-1} \in \{0,1\}$ we have $\{ q^{(i-1)}_{s_{i-1}}, q^{(i-1)}_{s'_{i-1}} \} = \{ q^{(i-1)}_0, q^{(i-1)}_1 \}$, thus by $\eqref{del15}$ (from Lemma $\ref{fo}$)
	\[ d( q^{(i-1)}_{s_{i-1}}, q^{(i-1)}_{s'_{i-1}} ) = 5 \delta_{i-1}. \]
	Now, since 
	\[ q^{(0)}_{s_{0}} q^{(1)}_{s_{1}} \dots q^{(i-2)}_{s_{i-2}} = q^{(0)}_{s'_{0}}q^{(1)}_{s'_{1}} \dots q^{(i-2)}_{s'_{i-2}}, \]
	(as $s\restrict{i-1} = s'\restrict{i-1}$) the invariance of $d$ implies 
	\begin{equation} \label{ctav}
	d\left( q^{(0)}_{s_{0}}q^{(1)}_{s_{1}} \dots q^{(i-2)}_{s_{i-2}} \cdot q^{(i-1)}_{s_{i-1}} , q^{(0)}_{s'_{0}} q^{(1)}_{s'_1} \dots q^{(i-2)}_{s'_{i-2}} \cdot q^{(i-1)}_{s'_{i-1}} \right) = 5\delta_{i-1},
	\end{equation}
	which proves $\eqref{kell}$, thus finishes the proof of $\eqref{diszj}$. Observe that (by the definition of the $C_s$'s $\eqref{Cs}$) $\eqref{ctav}$ implies the following
	\begin{equation} \label{Ctav}
	s\neq s' \in 2^i, s\restrict{i-1} = s'\restrict{i-1} \ (i>0) \ \rightarrow \ \  \ d(C_s,C_{s'}) \geq 3 \delta_{i-1}  \ \wedge \ \diam (C_s \cup C_{s'}) \leq 7 \delta_{i-1}.
	\end{equation}

	We define $\mu$ to be the standard coin-tossing measure supported by $C$, i.e. if $s \in 2^i$ then $\mu(C_s) = \frac{1}{2^{|s|}}$ (and $\mu(G \setminus C) = 0$). The following claim together with its corollary will
	ensure that each translate of the closed set $F = \varphi(c) \su G$ is null w.r.t. $\mu$, that is, $\mu$ witnesses that $F = \varphi(c)$ is Haar null.
	\bc \label{metszo}
		Let $h,h' \in G$ be fixed. Let $i \in \omega$, $i>0$. Then there exists $t \in 2^i$ such that
		\[ \left( h\cdot F \cdot h' \right) \cap \left(\bigcup_{s \in 2^i} C_s\right) \su C_t. \] 
	\ec
	Before proving the claim first observe that since $\mu(C_s) = \frac{1}{2^{|s|}}$ ($s \in 2^{<\omega}$)
	Claim $\ref{metszo}$ has the following corollary.
	\bcor
		Let $h,h' \in G$ be fixed, then
		\[ \mu(h \cdot F \cdot h') = 0, \]
		therefore $F = \varphi(c)$  is Haar null.
	\ecor
	\bp(Claim $\ref{metszo}$)
		Recall the definition of $\varphi(c) = F = \bigcap_{i \in \omega} \bigcup_{\psi \in \mathcal{T}_i(c)} F^\psi_c$ (Definition $\ref{fc}$).
		We will prove for this fixed $c \in 2^\omega$, $h,h' \in G$ that  for each $i \in \omega$
		\begin{equation} \label{indi}
		\exists t \in 2^{i}, \ \exists \pi \in \mathcal{T}_i(c): \ \  h\left( \bigcup_{\psi \in \mathcal{T}_i(c)} F^\psi_c \right)h' \cap C^i \su C_t \cap hF^\pi_ch'
		\end{equation}
		by induction on $i$.
		Before proving $\eqref{indi}$, we briefly describe the idea: In the $i$-th step we will use the induction hypothesis, that is $t' \in 2^{i-1}$, $\pi' \in \mathcal{T}_{i-1}(c)$ are such that
		$h\left( \bigcup_{\psi \in \mathcal{T}_{i-1}(c)} F^\psi_c \right)h' \cap C^{i-1} \su C_{t'} \cap hF^{\pi'}_ch'$
		and we will prove and use the fact that the $hF^\psi_ch'$'s ($\psi \in \mathcal{T}_i(c)$) are small enough (compared to $d(C_{t'\tieconcat 0}, C_{t'\tieconcat 1})$) so that each $hF^\psi_ch'$ can intersect only one of the $C_{t' \tieconcat j}$'s. Moreover, we will verify that the
		$hF^\psi_ch'$'s are far enough from each other so that whenever $\eta \in \mathcal{T}_i(c)$ is such that $hF^\eta_ch' \cap (C_{t'\tieconcat0} \cup C_{t' \tieconcat1}) \neq \emptyset$, then no $hF^\psi_ch'$ ($\psi \neq \eta$) can intersect $C_{t'\tieconcat0} \cup C_{t' \tieconcat1}$.
		
		Now turning to the proof of $\eqref{indi}$, for $i = 0$, let $t = \emptyset$ be the only element of $2^0$ (hence $C^0 = C_\emptyset$), and $\pi$ be the only element of $\mathcal{T}_0(c)$ (Remark $\ref{mcT_0(c)}$), 
		thus $\eqref{indi}$ obviously holds. 
		
		Suppose that $i>0$ and $t' \in 2^{i-1}$, $\pi' \in \mathcal{T}_{i-1}(c)$ are such that
		\[ h\left( \bigcup_{\psi \in \mathcal{T}_{i-1}(c)} F^\psi_c \right)h' \cap C^{i-1} \su C_{t'} \cap hF^{\pi'}_ch'. \]
		By Corollary $\ref{Ffa2}$ 
		\[ h\left( \bigcup_{\psi \in \mathcal{T}_{i}(c)} F^\psi_c\right) h' \su h\left( \bigcup_{\psi \in \mathcal{T}_{i-1}(c)} F^\psi_c \right) h',\] 
		implying
		\begin{equation} \label{ind00} h\left( \bigcup_{\psi \in \mathcal{T}_{i}(c)} F^\psi_c \right)h' \cap C^{i-1} \su C_{t'} \cap hF^{\pi'}_ch'. \end{equation}
		Using that the $C^i$'s are decreasing $\eqref{Cimon}$
		\begin{equation} \label{ind0} h\left( \bigcup_{\psi \in \mathcal{T}_{i}(c)} F^\psi_c \right)h' \cap C^{i} \su C_{t'} \cap hF^{\pi'}_ch'. \end{equation}
		Now using only that $h\left( \bigcup_{\psi \in \mathcal{T}_{i}(c)} F^\psi_c \right)h' \cap C^{i} \su C_{t'}$, after intersecting both sides with $C^i$ we obtain  (recalling that $C^i \cap C_{t'} = C_{t'\tieconcat 0} \cup C_{t'\tieconcat 1}$ by $\eqref{C^i}$ and $\eqref{diszj}$) the following:
		\begin{equation} \label{ind}  h\left( \bigcup_{\psi \in \mathcal{T}_{i}(c)} F^\psi_c \right)h' \cap C^{i} \su C^i \cap C_{t'} = C_{t'\tieconcat 0} \cup C_{t'\tieconcat 1}. \end{equation}	    
		Recall that Corollary $\ref{Ffa}$ implies that $(\bigcup_{\psi \in \mathcal{T}_{i}(c)} F^\psi_c) \cap F^{\pi'}_c = \bigcup_{\psi \in \mathcal{T}_i(c), \psi \geq \pi'} F^\psi_c$ from which
		$h(\bigcup_{\psi \in \mathcal{T}_{i}(c)} F^\psi_c)h' \cap hF^{\pi'}_ch' = h\left( \bigcup_{\psi \in \mathcal{T}_i(c), \psi \geq \pi'} F^\psi_c \right)h'$.
		These together with the inclusion $h\left( \bigcup_{\psi \in \mathcal{T}_{i}(c)} F^\psi_c \right)h' \cap C^{i} \su  hF^{\pi'}_ch'$ (that holds by $\eqref{ind0}$) imply
		\begin{equation} \label{uteq} \begin{array}{c} h\left( \bigcup_{\psi \in \mathcal{T}_{i}(c)} F^\psi_c \right)h' \cap C^{i} \su h\left( \bigcup_{\psi \in \mathcal{T}_{i}(c)} F^\psi_c \right)h' \cap hF^{\pi'}_ch' \su \\ \su h\left(\bigcup_{\psi \in \mathcal{T}_i(c), \psi \geq \pi'} F^\psi_c\right)h' = \left(\bigcup_{\psi \in \mathcal{T}_i(c), \psi \geq \pi'} hF^\psi_ch'\right). \end{array} \end{equation}
		We can summarize $\eqref{ind}$ and $\eqref{uteq}$ in
		\begin{equation} \label{ind+} h\left( \bigcup_{\psi \in \mathcal{T}_{i}(c)} F^\psi_c \right)h' \cap C^{i} \su (C_{t'\tieconcat 0} \cup C_{t'\tieconcat 1}) \cap \left(\bigcup_{\psi \in \mathcal{T}_i(c), \psi \geq \pi'} hF^\psi_ch'\right). \end{equation}
		
		Suppose that there is a $\psi \in \mathcal{T}_i(c), \psi \geq \pi'$ such that 
		$\left(C_{t'\tieconcat 0} \cup C_{t'\tieconcat 1}\right) \cap hF^\psi_ch' \neq \emptyset$, and let $\pi$  be such a $\psi$ if there is one,
		otherwise $h\left( \bigcup_{\psi \in \mathcal{T}_{i}(c)} F^\psi_c \right)h' \cap C^{i}$ is empty, (in this case let $\pi \in \mathcal{T}_i(c)$ be arbitrary, $\eqref{indi}$ obviously holds). This means that from now on we can assume that $\left(C_{t'\tieconcat 0} \cup C_{t'\tieconcat 1}\right) \cap hF^\pi_ch' \neq \emptyset$,
		let $t \in \{ t'\tieconcat 0, t' \tieconcat 1\}$ be such that $hF^\pi_ch' \cap C_{t} \neq \emptyset$.
		
		It remains to show that 
		\begin{equation} \label{Hnullutso} \left(C_{t'\tieconcat 0} \cup C_{t'\tieconcat 1}\right) \cap \left(\bigcup_{\psi \in \mathcal{T}_i(c), \psi \geq \pi'} hF^\psi_ch' \right)  \su hF_c^\pi h' \cap C_t \end{equation} 
		(this with $\eqref{ind+}$ will complete the proof of $\eqref{indi}$).
		Corollary $\ref{Ffa}$ gives that for  $\psi \in \mathcal{T}_i(c)$ $\psi \geq \pi'$, if $\psi \neq \pi$ then $d(F_c^\psi, F_c^\pi) > 7 \delta_{i-1}$ (that is equal to $d(hF_c^\psi h', hF_c^\pi h') > 7 \delta_{i-1}$ by the invariance of $d$). By $\eqref{Ctav}$ $\diam(C_{t' \tieconcat 0} \cup C_{t' \tieconcat 1}) \leq 7 \delta_{i-1}$, therefore
		\begin{equation} \label{nemmetsz} \psi \in \mathcal{T}_i(c), \psi \geq \pi', \  \psi \neq \pi \ \ \to hF_c^\psi h' \cap (C_{t' \tieconcat 0} \cup C_{t' \tieconcat 1})) = \emptyset. \end{equation}
		On the other hand, since $\pi \in \mathcal{T}_i(c)$, $\diam(hF_c^\pi h') \leq 2 \delta_{i-1}$ by Definition $\ref{Fdef}$ and the invariance of $d$, and $d(C_{t' \tieconcat 0},C_{t' \tieconcat 1}) \geq 3\delta_{i-1}$ by $\eqref{Ctav}$, therefore $hF_c^\psi h'$ can only intersect one of $C_{t' \tieconcat 0}$, $C_{t' \tieconcat 1}$. Thus $hF^\pi_ch' \cap C_{t} \neq \emptyset$ together with $\eqref{nemmetsz}$ verifies $\eqref{Hnullutso}$.
	\ep

Next we prove that perfectly many $\varphi(c)$ form a right compact catcher set.
\bl \label{lemm}
	Let $P \su 2^\omega$ be non-empty, perfect, and $K \su G$ be compact. Then there exists $h \in G$ such that
	\begin{equation} \label{jobbtolt} Kh \su \bigcup_{c \in P} \varphi(c). \end{equation}
\el
\bp
	We will need the following technical statements.
	Let $T \su 2^{<\omega}$ denote the downward closed tree corresponding to $P$, i.e.
	\begin{equation} \label{Pdef} P = [T] = \{ c \in 2^\omega: \ \ \forall i\in \omega \ c\restrict{i} \in T \} \end{equation}
	(see \cite[(2.4)]{Ke}).
	
	We will need the following Lemma.
	\bl \label{vegtelensorozat}
		There exist sequences $\underline{n} = (n_j)_{j \in \omega}, \underline{m} = (m_j)_{j \in \omega} \in \omega^\omega$, $(\psi_i)_{i \in \omega}$ such that the following hold.
		\begin{enumerate}[(i)]
			\item \label{k1} $\psi_i \in \mathcal{T}_i$,
			\item \label{alphaii+1} $\psi_{i-1} \leq \psi_{i}$ ($i>0$),
			\item \label{T'_i} $\dom( \psi_i ) = \cup_{l\leq i} \prod_{j<l} n_j$, (i.e. $\underline{n}^{\psi_i}\restrict{i} = \underline{n}\restrict{i}$),
			\item \label{reszeT} $\ran (\psi_i) \su T$,
			\item \label{kut-1} $m^{\psi_i}_j = m_j$ ($j<i$),
			\item \label{kut} for each $i \in \omega$, if $i>0$ then
			\begin{equation} \begin{array}{l}
			Kt^{(0)}_{n_{0},m_{0}} \cdot t^{(1)}_{n_{1},m_{1}} \cdot \dots \cdot  t^{(i-1)}_{n_{i-1},m_{i-1}} 
			\su  \bigcup_{s \in \prod_{j<i} n_j} B(\frac{\delta_{i-1}}{2}) g_{\underline{n}\restrict{i}, \underline{m}\restrict{i}, s} = \\
			= \bigcup_{s \in \prod_{j<i} n_j} B(\frac{\delta_{i-1}}{2}) \left(  q^{(0)}_{s_{0}} t^{(0)}_{n_{0},m_{0}} \right) \cdot \left( q^{(1)}_{s_{1}} t^{(1)}_{n_{1},m_{1}} \right) \cdot \dots  \cdot \left( q^{(i-1)}_{s_{i-1}} t^{(i-1)}_{n_{i-1},m_{i-1}} \right) .
			\end{array}
			\end{equation}
				\end{enumerate}
	\el
	\bp
		For $i=0$ let $\psi_0: \{ \emptyset \} \to \{ \emptyset \}$ be the unique element of $\mathcal{T}_0$. It is straightforward to check that $\eqref{k1}-\eqref{kut}$  hold.
		
		Suppose that $i>0$ and $\psi_j$ ($j<i$) and $m_j, n_j$ ($j<i-1$) are already constructed satisfying $\eqref{k1}-\eqref{kut}$.
		By $\eqref{kut}$
		\begin{equation} \label{kuti-1}   \begin{array}{l}
		Kt^{(0)}_{n_{0},m_{0}} \cdot t^{(1)}_{n_{1},m_{1}} \cdot \dots \cdot  t^{(i-2)}_{n_{i-2},m_{i-2}}  \su 
		\bigcup_{s \in \prod_{j<i-1} n_j} B(\frac{\delta_{i-2}}{2}) \cdot g_{\underline{n}\restrict{i-1}, \underline{m}\restrict{i-1}, s}.
		\end{array} \end{equation}
		The definition of the $q^{(i-1)}_j$'s (Lemma $\ref{fo}$ $\eqref{gk}$) implies that $\{ q^{(i-1)}_j: \ j \in \omega \} \su B(\delta_{i-2}/2)$ is dense in $B(\delta_{i-2}/2)$, therefore  (for $s \in \prod_{j<i-1} n_j$)
		\[  g_{\underline{n}\restrict{i-1}, \underline{m}\restrict{i-1}, s} \cdot \{ q^{(i-1)}_j: \ j \in \omega \} \]
		\[ \text{ is dense in } \]
		\[    g_{\underline{n}\restrict{i-1}, \underline{m}\restrict{i-1}, s} \cdot B(\delta_{i-2}/2 ) =  B(\delta_{i-2}/2) \cdot   g_{\underline{n}\restrict{i-1}, \underline{m}\restrict{i-1}, s}\]
		(where the last equality is due to Lemma $\ref{gomb}$).
		Now because translations of a fixed open set by a dense set cover we have that (for a fixed $s \in \prod_{j<i-1} n_j$)
		\[ B(\delta_{i-2}/2) \cdot g_{\underline{n}\restrict{i-1}, \underline{m}\restrict{i-1}, s} \su \bigcup_{j \in \omega} B(\delta_{i-1}/2) \cdot g_{\underline{n}\restrict{i-1}, \underline{m}\restrict{i-1}, s} \cdot q_j^{(i-1)}.  \]

		Therefore by $\eqref{kuti-1}$ and the compactness of $K$, there exist $l \in \omega$ such
		that
		\begin{equation} \label{Kkozb} \begin{array}{l} 	Kt^{(0)}_{n_{0},m_{0}} \cdot t^{(1)}_{n_{1},m_{1}} \cdot \dots \cdot  t^{(i-2)}_{n_{i-2},m_{i-2}}
		\su \bigcup_{ s \in \prod_{j<i-1} n_j, k<l}  B(\delta_{i-1}/2) \cdot g_{\underline{n}\restrict{i-1}, \underline{m}\restrict{i-1}, s} \cdot q^{(i-1)}_k. \end{array} \end{equation}
		Define $n_{i-1} = l$.
		
		Now we are ready for the construction of $\psi_i$.
		For each $s \in \prod_{j<i-1} n^{\psi_{i-1}}_j = \prod_{j<i-1} n_j $ choose pairwise incomparable proper extensions $s_0,s_1, \dots, s_{n_{i-1}-1}$ of $ \psi_{i-1}(s)$ in $T$, i.e.
		\begin{equation} \label{tk}
		\left( \begin{array}{l}  \psi_{i-1}(s) \subsetneq s_j \in T, \text{ and} \\
		 \text{if }k \neq j < n_{i-1} \text{, then} \ \ s_k \nsubseteq s_j \ \wedge \ s_j \nsubseteq s_k \end{array} \right) \  \left(\forall s \in \prod_{j<i-1} n_j\right).  \end{equation}
		(This can be done since $\ran(\psi_{i-1})\su T$, and $T$ is a perfect tree, i.e. each of its element has incomparable extensions, see \cite[(6.14)]{Ke}.)
		Define 
		\begin{equation} \label{psibov}
		\begin{array}{l} \psi_i \supset  \psi_{i-1}, \ \ \dom(\psi_{i}) = \bigcup_{l\leq i} \prod_{j<l} n_j = \dom(\psi_{i-1}) \cup \prod_{j<i} n_j, \\ 
		\psi_i(s \tieconcat k) =  s_k \ \ (k<n_{i-1}).	\end{array} \end{equation}
		Now we are ready to check $\eqref{k1}-\eqref{reszeT}$ for $\psi_i$ and the $n_j$'s ($j<i$).
		$\eqref{T'_i}$ obviously holds, and it is straightforward to check $\eqref{k1}$, $\eqref{alphaii+1}$ (recalling the definition of $\mathcal{T}$ and the partial order on it, Definitions $\ref{Tdef}$ and $\ref{Tpo}$).
		Using $\eqref{reszeT}$ for $\psi_{i-1}$, and that the $s_k$'s ($s \in \prod_{j<i-1} n_j$) are in $T$ $\eqref{tk}$, we have that $\ran(\psi_i) \su T$  by $\eqref{psibov}$, i.e. $\eqref{reszeT}$ holds for $\psi_i$. 

		Now, as $\psi_{i} \in \mathcal{T}_i$, put $m_{i-1} = m^{\psi_{i}}_{i-1}$. By the definition of the $\underline{m}^\psi$'s (Definition $\ref{mdef}$) $\psi_{i-1} \leq \psi_i$ implies that $(\underline{m}^{\psi_{i}})\restrict{i-1} = \underline{m}^{\psi_{i-1}}$. But $m^{\psi_{i-1}}_j =m_j$ for $j<i-1$ by the induction hypothesis $\eqref{kut-1}$, thus $\eqref{kut-1}$ holds with $\psi_i$, too.
		
		It remains to check $\eqref{kut}$. As we have already defined $n_{i-1}$ to be such that $l = n_{i-1}$ satisfies $\eqref{Kkozb}$, we can formulate it as
		\begin{equation} \label{Kkoz2}	Kt^{(0)}_{n_{0},m_{0}} \cdot t^{(1)}_{n_{1},m_{1}} \cdot \dots \cdot  t^{(i-2)}_{n_{i-2},m_{i-2}}
		\su \bigcup_{ s \in \prod_{j<i} n_j}  B(\delta_{i-1}/2) \cdot g_{\underline{n}\restrict{i-1}, \underline{m}\restrict{i-1}, s\restrict{i-1}} \cdot q^{(i-1)}_{s_{i-1}}. \end{equation}
		Using that for $s \in \prod_{j<i} n_j$  (by Definition $\ref{gdef}$)
		\[ g_{\underline{n}\restrict{i-1}, \underline{m}\restrict{i-1}, s\restrict{i-1}} \cdot q^{(i-1)}_{s_{i-1}}\cdot t^{(i-1)}_{n_{i-1},m_{i-1}} = g_{\underline{n}\restrict{i}, \underline{m}\restrict{i}, s}, \]
		the multiplication of $\eqref{Kkoz2}$ by $t^{(i-1)}_{n_{i-1},m_{i-1}}$ from the right gives
		\[ Kt^{(0)}_{n_{0},m_{0}} \cdot t^{(1)}_{n_{1},m_{1}} \cdot \dots \cdot  t^{(i-2)}_{n_{i-2},m_{i-2}} t^{(i-1)}_{n_{i-1},m_{i-1}} 
		\su \bigcup_{ s \in \prod_{j<i} n_j}  B(\delta_{i-1}/2) \cdot g_{\underline{n}\restrict{i}, \underline{m}\restrict{i}, s}, \]
		as desired.
	\ep

	For proving $\eqref{jobbtolt}$ our candidate for $h$ is 
	$\lim_{i \to \infty} \left( t^{(0)}_{n_{0},m_{0}} \cdot t^{(1)}_{n_{1},m_{1}} \cdot \dots \cdot  t^{(i-1)}_{n_{i-1},m_{i-1}} \right)$ given by Lemma $\ref{vegtelensorozat}$.
	First we have to prove that this limit exists.

	\bl \label{sortav}
		For $i<k \in \omega$ and any sequences $0 < p_i, p_{i+1}, \dots, p_k$, $0 \leq r_i, r_{i+1}, \dots, r_k$
		\[ t^{(i)}_{p_i,r_i} t^{(i+1)}_{p_{i+1},r_{i+1}} \ldots t^{(k)}_{p_k,r_k} \in B(\delta_{i-1}/2). \]
	\el
	\bp
		By Lemma $\ref{eltolos}$ $t^{(j)}_{p_j,r_j} \in B(\varepsilon_j)$, hence Lemma $\ref{szorzasfact}$ implies that
		\begin{equation} \label{eps0} d (t^{(i)}_{p_i,r_i} t^{(i+1)}_{p_{i+1},r_{i+1}} \ldots t^{(k)}_{p_k,r_k},e) < \varepsilon_i + \varepsilon_{i+1} +\ldots +  \varepsilon_{k}, \end{equation}
		therefore we only have to prove that
		\begin{equation} \label{eps} \varepsilon_i + \varepsilon_{i+1} + \ldots + \varepsilon_{k} \leq \delta_{i-1}/2. \end{equation}
		But $\varepsilon_l + \delta_l < \delta_{l-1}/2$ by $\eqref{del1}$ from Lemma $\ref{fo}$, from which 
		(by induction for $j \geq i$, always replacing $\delta_j$ by the smaller $\varepsilon_{j+1} + \delta_{j+1}$)
		\[ \varepsilon_i + \varepsilon_{i+1} + \ldots + \varepsilon_{j} + \delta_j  < \delta_{i-1}/2. \]
		From this $\eqref{eps}$ follows.
			\ep
	
	\bcor
		
		For the sequence $\left( t^{(0)}_{n_{0},m_{0}} \cdot t^{(1)}_{n_{1},m_{1}} \cdot \ldots \cdot  t^{(i-1)}_{n_{i-1},m_{i-1}} \right)_{i \in \omega}$ (given by Lemma $\ref{vegtelensorozat}$), if $i \leq j$
		\[ d(t^{(0)}_{n_{0},m_{0}} \cdot t^{(1)}_{n_{1},m_{1}} \cdot \ldots \cdot  t^{(i-1)}_{n_{i-1},m_{i-1}}, t^{(0)}_{n_{0},m_{0}} \cdot t^{(1)}_{n_{1},m_{1}} \cdot \ldots \cdot  t^{(j-1)}_{n_{j-1},m_{j-1}})  < \delta_{i-1}/2, \]
		in particular, it is Cauchy.
	\ecor
	\bp
		By the invariance of $d$
		\[ d(t^{(0)}_{n_{0},m_{0}} \cdot t^{(1)}_{n_{1},m_{1}} \cdot \ldots \cdot  t^{(i-1)}_{n_{i-1},m_{i-1}}, t^{(0)}_{n_{0},m_{0}} \cdot t^{(1)}_{n_{1},m_{1}} \cdot \ldots \cdot  t^{(j-1)}_{n_{j-1},m_{j-1}})  = \]
		\[ d(e, t^{(i)}_{n_{i},m_{i}} \cdot t^{(i+1)}_{n_{i+1},m_{i+1}} \cdot \ldots \cdot  t^{(j-1)}_{n_{j-1},m_{j-1}}), \]
		therefore Claim $\ref{sortav}$ yields the inequality (furthermore, as the sequence $(\delta_i)_{i\in \omega}$ tends to $0$ it is Cauchy, indeed).
	\ep
	
	Now we can define $h$.
	\begin{equation} \label{hdef}
	h = \lim_{i \to \infty} \left( t^{(0)}_{n_{0},m_{0}} \cdot t^{(1)}_{n_{1},m_{1}} \cdot \dots \cdot  t^{(i-1)}_{n_{i-1},m_{i-1}} \right).
	\end{equation}
	
	The following two lemmas will ensure that $\eqref{jobbtolt}$ holds.
	\bl \label{elsotart}
		With $(n_i)_{i \in \omega}$, $(m_i)_{i \in \omega}$ given by Lemma $\ref{vegtelensorozat}$
		and $h$ as in $\eqref{hdef}$,
		\[ \begin{array}{l} Kh \su \bigcap_{i>0}  \bigcup_{s \in \prod_{j<i} n_j} \cl(B(\delta_{i-1})) g_{\underline{n}\restrict{i}, \underline{m}\restrict{i}, s} = \\
		=  \bigcap_{i >0} \bigcup_{s \in \prod_{j<i} n_j} \cl(B(\delta_{i-1})) \left(  q^{(0)}_{s_{0}} t^{(0)}_{n_{0},m_{0}} \right) \cdot \left( q^{(1)}_{s_{1}} t^{(1)}_{n_{1},m_{1}} \right) \cdot \dots  \cdot \left( q^{(i-1)}_{s_{i-1}} t^{(i-1)}_{n_{i-1},m_{i-1}} \right). 
		
		\end{array} \]
		
	\el
	\bl \label{masodiktart}
		\[ \bigcap_{i >0} \bigcup_{s \in \prod_{j<i} n_j} \cl(B(\delta_{i-1})) g_{\underline{n}\restrict{i}, \underline{m}\restrict{i}, s} \su \bigcup_{c \in P} \varphi(c). \]
	\el
	\bp(\textnormal{Lemma} $\ref{elsotart}$)
		Fix $i \in \omega$, $i > 0$.
		By $\eqref{kut}$ we know that
		\begin{equation} \label{kut'} \begin{array}{l}
		Kt^{(0)}_{n_{0},m_{0}} \cdot t^{(1)}_{n_{1},m_{1}} \cdot \dots \cdot  t^{(i-1)}_{n_{i-1},m_{i-1}} 
		\su  \bigcup_{s \in \prod_{j<i} n_j} B(\frac{\delta_{i-1}}{2}) g_{\underline{n}\restrict{i}, \underline{m}\restrict{i}, s} = \\
		= \bigcup_{s \in \prod_{j<i} n_j} B(\frac{\delta_{i-1}}{2}) \left( q^{(0)}_{s_{0}} t^{(0)}_{n_{0},m_{0}}  \right) \cdot \left( q^{(1)}_{s_{1}} t^{(1)}_{n_{1},m_{1}}\right) \cdot \dots  \cdot \left( q^{(i-1)}_{s_{i-1}} t^{(i-1)}_{n_{i-1},m_{i-1}} \right) .
		\end{array}
		\end{equation}
		Now Claim $\ref{sortav}$ states that if $k\geq i$, then
		\[ 
		t^{(i)}_{n_i,m_i} t^{(i+1)}_{n_{i+1},m_{i+1}} \ldots t^{(k)}_{n_k,m_k} \in B(\delta_{i-1}/2),
		\]
		therefore 
		\begin{equation} \label{h'} h' = \lim_{l \to \infty} t^{(i)}_{n_i,n_i} t^{(i+1)}_{n_{i+1},m_{i+1}} \ldots t^{(i+l)}_{n_{i+l},m_{i+l}} \in \cl(B(\delta_{i-1}/2)). \end{equation}
		Observe that 
		\begin{equation} \label{hh'} \begin{array}{c} h 
		= t^{(0)}_{n_{0},m_{0}} \cdot t^{(1)}_{n_{1},m_{1}} \cdot \dots \cdot  t^{(i-1)}_{n_{i-1},m_{i-1}} \cdot h'. \end{array} \end{equation}
		
		Using $\eqref{hh'}$, $\eqref{kut'}$ and finally  $\eqref{h'}$ 
		\[  Kh =	Kt^{(0)}_{n_{0},m_{0}} \cdot t^{(1)}_{n_{1},m_{1}} \cdot \dots \cdot  t^{(i-1)}_{n_{i-1},m_{i-1}} \cdot h'
		\su  \left( \bigcup_{s \in \prod_{j<i} n_j} B(\frac{\delta_{i-1}}{2}) g_{\underline{n}\restrict{i}, \underline{m}\restrict{i}, s} \right) h'  \su \]
		\[ \su  \bigcup_{s \in \prod_{j<i} n_j} \left (B(\delta_{i-1}/2) \cdot g_{\underline{n}\restrict{i}, \underline{m}\restrict{i}, s} \cdot \cl(B(\delta_{i-1}/2)) \right) \su \]
		\[ \su \bigcup_{s \in \prod_{j<i} n_j} \cl \left( B(\delta_{i-1}/2) \cdot g_{\underline{n}\restrict{i}, \underline{m}\restrict{i}, s} \cdot B(\delta_{i-1}/2) \right) \]
		and (since by Lemmas $\ref{gomb}$, $\ref{szorzasfact}$ 
		$$ B(\delta_{i-1}/2) \cdot g \cdot B(\delta_{i-1}/2) = B(\delta_{i-1}/2) \cdot B(\delta_{i-1}/2) g 
		\su B(\delta_{i-1}) g \ (g \in G)$$
		holds,)
		\[ Kh  \su	\bigcup_{s \in \prod_{j<i} n_j} \cl \left( B(\delta_{i-1}) \cdot g_{\underline{n}\restrict{i}, \underline{m}\restrict{i}, s} \right) = \bigcup_{s \in \prod_{j<i} n_j} \cl \left( B(\delta_{i-1}) \right) \cdot g_{\underline{n}\restrict{i}, \underline{m}\restrict{i}, s}, \]
		we are done.
	\ep
	
	\bp(\textnormal{Lemma} $\ref{masodiktart}$)
		Recall that $T \su 2^{<\omega}$ is defined so that the perfect set $P \su 2^{\omega}$ is the body of $T$ ($\eqref{Pdef}$),
		and the sequences $(n_i)_{i \in \omega}, (m_i)_{i \in \omega} \in \omega^\omega$, $(\psi)_{i \in \omega}$ fulfill the criteria $\eqref{k1}-\eqref{kut}$.
		
		Fix an element $x \in \bigcap_{i >0} \bigcup_{s \in \prod_{j<i} n_j} \cl(B(\delta_{i-1})) g_{\underline{n}\restrict{i}, \underline{m}\restrict{i}, s}$,
		and define $B \su \bigcup_{i>0} \prod_{j<i}n_j$ such that
		\[ B = \{r \in \bigcup_{i>0} \prod_{j<i}n_j: \  x \in \cl(B(\delta_{i-1})) g_{\underline{n}\restrict{i}, \underline{m}\restrict{i}, r}\}.\] 
		First observe that using Lemma $\ref{monotonGomb}$ $B$ is downward closed. $B$ is an infinite tree by the choice of $x$, and it has finitely many nodes on each level, thus by K\H onig's Lemma
		there is an infinite branch through it. Let $r' \in \prod_{i<\omega} n_i$ be the union of that branch, i.e. the corresponding infinite sequence, for which
		\begin{equation} \label{r'} x \in \cl(B(\delta_{i-1})) g_{\underline{n}\restrict{i}, \underline{m}\restrict{i}, r'\restrict{i}} \ \ (\forall i >0). \end{equation}
		Now, for $i \leq k$, $  \psi_i \leq  \psi_k$ by $\eqref{alphaii+1}$ from Lemma $\ref{vegtelensorozat}$, that implies (using Definition $\ref{Tpo}$) that
		$\psi_i = (\psi_k)\restrict{\dom(\psi_i)}$. Recall that $\dom(\psi_k) = \bigcup_{l\leq k} \prod_{j<l} n_j$ by $\eqref{T'_i}$ from Lemma $\ref{vegtelensorozat}$.
		This implies that for $i < k$
		\[  \psi_i(r'\restrict{i}) = \psi_k(r'\restrict{i}) \subsetneq \psi_k(r'\restrict{k}). \]
		Define 
		\begin{equation} \label{c'def} c' = \bigcup \{\psi_i(r'\restrict{i}): \ i \in \omega \}. \end{equation}
		Then $c' \in P$, since $\ran(\psi_i) \su T$ ($i \in \omega$) by $\eqref{reszeT}$ from Lemma $\ref{vegtelensorozat}$ and  $\eqref{Pdef}$.
		Now we will verify that $x \in \varphi(c')$.
		Recall the definition of $\varphi(c')$ (Definition $\ref{varphidef}$):
		\[ \begin{array}{l} \varphi(c') = \bigcap_{i \in \omega} \bigcup_{\psi \in \mathcal{T}_i(c')} F^\psi_{c'} = \bigcap_{i >0} \bigcup_{\psi \in \mathcal{T}_i(c')} \cl(B(\delta_{i-1})) g_{\underline{n}^\psi,\underline{m}^\psi,b^{c',\psi}}. \\
		\end{array} \]
		Hence for proving that $x \in \varphi(c')$ it suffices (by $\eqref{r'}$) to show that for each $i \in \omega$,  if $i>0$ then $\psi_i \in \mathcal{T}_i(c')$ and $b^{c',\psi_i} = r'\restrict{i}$
		(since we know that $\underline{n}^{\psi_i} = \underline{n}\restrict{i}$, $\underline{m}^{\psi_i} = \underline{m}\restrict{i}$ by $\eqref{kut-1}$ and $\eqref{T'_i}$ from Lemma $\ref{vegtelensorozat}$).
		
		Fix $i \in \omega$, $i>0$. Now $\psi_i(r'\restrict{i}) \subset c'$ by the definition of $c'$ $\eqref{c'def}$, therefore $\psi_i \in \mathcal{T}_i(c')$ by Definition $\ref{cfaj}$, and $r'\restrict{i} = b^{c',\psi_i}$ (by Definition $\ref{ccdef}$).
		\ep

\bl
	The mapping $\varphi: 2^\omega \to \mathcal{F}(G)$ is Borel.
\el
We need that for each Borel set $B \su \mathcal{F}(G)$ (w.r.t. the Effros Borel structure) its preimage, $\varphi^{-1}(B) \in \mathcal{B}(2^\omega)$. For this it suffices to show that for each fixed open set $U \su G$ the preimage of
\[ \{ F \in \mathcal{F}(G): \ F \cap U \neq \emptyset \} \]
under $\varphi$ is a Borel subset of $2^\omega$, since those sets form a generator system of the Effros Borel structure on $\mathcal{F}(G)$. This is provided by the following Lemma.
\bl \label{Borelcl}
	Let the open set $U \su G$ be fixed. Suppose that  $c \in 2^\omega$ is such that $U \cap \varphi(c) \neq \emptyset$.
	Then there exists $i_0 \in \omega$ such that for each $c' \in [c\restrict{i_0}] = \{ y \in 2^\omega: \ y\restrict{i_0} = c\restrict{i_0} \}$
	we have $\varphi(c') \cap U \neq \emptyset$. 
	In particular 
	\[ \varphi^{-1}( \{F \in \mathcal{F}(G): \ F \cap U  \neq \emptyset \}) \su 2^\omega \text{ is open.} \] 
\el
\bp
	First pick an element $x \in U \cap \varphi(c)$. By the definition of $\varphi(c)$ (Definition $\ref{varphidef}$) for each $i>0$
	\[ x \in \bigcup_{\psi \in \mathcal{T}_i(c)} F^\psi_c 
	.\]
	Now by Corollary $\ref{Ffa}$ for a fixed $i$ the $F^\psi_c$'s ($\psi \in \mathcal{T}_i(c)$) are disjoint, hence for each $i>0$ there is a unique 
	$\psi \in \mathcal{T}_i(c)$ for which 
	\[ x \in F^{\psi}_c 
	. \]
	As $\diam(F_c^\psi) \leq 2 \delta_{i-1}$ for $\psi \in \mathcal{T}_i(c)$ and $\delta_i$
	tends to $0$ as $i$ tends to $\infty$ ($\eqref{del1}$ from Lemma $\ref{fo}$), we can fix an index $i'$ (and the unique $\pi \in \mathcal{T}_{i'}(c)$) so that
	\begin{equation} \label{FreszeU} x \in F^{\pi}_c \su U. \end{equation}
	Observe that if $l$ witnesses that $\pi \in \mathcal{T}_{i'}(c)$ according to Definition $\ref{cfaj}$, i.e. 
	$c\restrict{l} = \pi(b^{c,\pi})$ (Definition $\ref{ccdef}$), then also for each  $y \in 2^\omega$ with $y\restrict{l} = c\restrict{l}$ 
	we have $\pi \in \mathcal{T}_i'(y)$, and $b^{y,\pi} = b^{c,\pi}$. Define $i_0 = l$.
	Fixing such a $y$ it follows from this that 
	\begin{equation} \label{Fegyenlo} F^\pi_c = F^\pi_y \end{equation} (Definition $\ref{Fdef}$). 
	
	Now we show that for each $y \in 2^\omega$ the tree $\mathcal{T}(y) \su \mathcal{T}$ is pruned, i.e. for each $i \in \omega$, $\psi \in \mathcal{T}_i(y)$ there exists
	$\psi' \in \mathcal{T}_{i+1}(y)$ with $\psi \leq \psi'$ (w.r.t. Definition $\ref{Tpo}$).  First we show that this  completes the proof of Claim $\ref{Borelcl}$.
	Indeed, an infinite branch 
	\[ \pi^{(0)} = \pi \leq \pi^{(1)} \leq \pi^{(2)} \leq \ldots \leq \pi^{(k)} \leq \ldots  \ \ (\pi^{(j)} \in \mathcal{T}_{i'+j}(y)) \]
	in $\mathcal{T}(y)$	would yield an infinite decreasing chain of closed sets
	\[ F^{\pi^{(0)}}_y = F^\pi_y \supset F_y^{\pi^{(1)}} \supset \ldots \supset  F_y^{\pi^{(k)}} \supset \ldots \]
	(using Corollary $\ref{Ffa2}$).
	Now $\pi^{(j)} \in \mathcal{T}_{i'+j}(y)$, $\diam(F_y^{\pi^{(j)}}) \leq 2 \delta_{i'+j-1}$ (by Definition $\ref{Fdef}$), and 
	the sequence $\delta_j$ ($j \in \omega$) tends to $0$ ($\eqref{del1}$ from Lemma $\ref{fo}$). 
	This would yield (recalling Definition $\ref{varphidef}$) that
	\[ \emptyset \neq \bigcap_{j \in \omega} F^{\pi^{(j)}}_y \su \varphi(y) = \bigcap_{i \in \omega} \bigcup_{\psi \in \mathcal{T}_i(y)} F^\psi_y, \]
	implying
	\[ \varphi(y) \cap F^\pi_y \neq \emptyset. \]
	This together with $\eqref{Fegyenlo}$ and $\eqref{FreszeU}$ implies
	\[ \varphi(y) \cap U \neq \emptyset, \]
	hence
	\[ \varphi^{-1}( \{F \in \mathcal{F}(G): \ F \cap U  \neq \emptyset \}) \supset \{ y \in 2^\omega: \ y\restrict{i_0} = c\restrict{i_0} \},  \]
	as desired.
	
	Finally we have to check that $\mathcal{T}(y)$ is pruned. Fix $i \in \omega$, $\psi \in \mathcal{T}_i(y)$. 
	According to Definition $\ref{ccdef}$ this means that in the set of maximal elements of $\dom(\psi)$ 
	\[ M = \prod_{j<i} n^\psi_j \]
	there is an element $s_0$ with $\psi(s_0) \su y$. 
	For each $s \in M$ choose $t_s \in 2^{<\omega}$ with $\psi(s) \subsetneq t_s$, and with $t_{s_0} \su y$.
	Define $\psi' \supset \psi$ so that 
	\[ \dom(\psi') = \left( \bigcup_{l \leq i} \prod_{j<l} n^\psi_j \right) \cup \left(\left(\prod_{j < i} n^\psi_j\right) \times 1\right)   \]
	(i.e $\underline{n}^{\psi'} = \underline{n}^\psi \tieconcat 1$, $\dom(\psi') =  \bigcup_{l \leq i+1} \prod_{j<l} n^{\psi'}_j$).
	For each $s \in M$ let  $\psi'(s \tieconcat 0) = t_{s}$. It is straightforward to check that $\psi' \in \mathcal{T}$ (Definition $\ref{Tdef}$), hence obviously $\psi' \in \mathcal{T}_{i+1}$ (Definition $\ref{T_k}$), and
	$\psi' \geq \psi$ (Definition $\ref{Tpo}$). Last, $\psi'(s_0 \tieconcat 0) \su y$ proves that $\psi' \in \mathcal{T}(y)$ (Definition $\ref{ccdef}$).
	\ep

\section{Open problems}

\bpr
Which results can one generalize to all Polish groups admitting a two-sided invariant metric?
\epr

\bpr
What can we say about the case of arbitrary Polish groups?
\epr

\section*{Acknowledgments} We are grateful for J. Brendle and J. Stepr\=ans for some illuminating discussions.

\end{document}